\documentclass[a4paper,10pt]{article}

\usepackage[utf8x]{inputenc}
\usepackage[T1]{fontenc} 
\usepackage{amsmath,amsthm}
\usepackage[dvips]{graphicx}
\usepackage{amsfonts,amssymb}
\usepackage{geometry}
\usepackage{hyperref}
\usepackage{stmaryrd}
\usepackage{fancyhdr}
\usepackage{url}

\newtheorem{prop}{Proposition}[section]
\newtheorem{LM}{Lemma}[section]
\newtheorem{thm}{Theorem}[section]

\newtheorem{df-prop}{Definition-Proposition}[section]
\newtheorem{cor}{Corollary}[section]

\newtheorem*{conj}{Conjecture}
\newtheorem*{claim}{Claim}
\newtheorem*{thm*}{Theorem \ref{distgen}}
\newtheorem*{thm**}{Theorem \ref{eqasai}}
\newtheorem*{thme**}{Theorem \ref{asaiprinc}}

\newtheoremstyle{pourlesremarques}{\topsep}{\topsep}{\normalfont}{}{\bfseries}{.}{ }{}
\theoremstyle{pourlesremarques}

\def\adots{\mathinner{\mkern2mu\raise 1pt\hbox{.}\mkern 3mu\raise
4pt\hbox{.}\mkern1mu\raise 7pt\hbox{{.}}}}

\title {\textbf{Distinction of some induced representations}}
\author{Nadir MATRINGE}

\begin{document}
 
\maketitle

\begin{abstract}
Let $K/F$ be a quadratic extension of $p$-adic fields, $\sigma$ the nontrivial element of the Galois group of $K$ over $F$, and $\Delta$ a 
quasi-square-integrable representation of $GL(n,K)$. Denoting by $\Delta^{\vee}$ the smooth contragredient of $\Delta$, and by $\Delta^{\sigma}$ the 
representation $\Delta\circ \sigma$, we show that representation of $GL(2n, K)$ obtained by normalized parabolic induction of the representation 
$\Delta^\vee \otimes \Delta^\sigma$, is distinguished with respect to $GL(2n,F)$. This is a step towards 
the classification of distinguished generic representations of general linear groups over $p$-adic fields.

\end{abstract}

\section*{Introduction}

Let $K/F$ be a quadratic extension of $p$-adic fields, $\sigma$ the nontrivial element of the Galois group of $K$ over $F$, and $\Delta$ a quasi-square-integrable 
representation of $GL(n,K)$. We denote by $\sigma$ again the automorphism of $M_{2n}(K)$ induced by $\sigma$.\\ 
If $\chi$ is a character of $F^*$, a smooth representation $\rho$ of $GL(2n,K)$ is said to be $\chi$-distinguished if there is a nonzero linear form $L$ on its 
space $V$, verifying $L(\rho(h)v)=\chi(det(h))L(v)$ for all $h$ in $GL(2n,F)$ and $v$ in $V$, we say distinguished if $\chi=1$. If $\rho$ is irreducible, the space of 
such linear forms is of dimension at most $1$ (Proposition 11 of \cite{F2}).\\
Calling $\Delta^\vee$ the smooth contragredient of $\Delta$ and $\Delta^\sigma$ the representation $\Delta\circ \sigma$, we denote by $\Delta^\sigma \times \Delta^\vee$ the 
representation of $GL(2n,K)$, obtained by normalized induction of the representation $\Delta^\sigma \otimes \Delta^\vee$ of the standard parabolic subgroup of type $(n,n)$.
 The aim of the present work is to show that the representation $\Delta^\sigma \times \Delta^\vee$ is distinguished.\\ 
The case $n=1$ is treated in $\cite{H}$ for unitary $\Delta^\sigma \times \Delta^\vee$, using a criterion characterizing distinction in terms of gamma factors. In 
\cite{F3}, Flicker defines a linear form on the space of $\Delta^\sigma \times \Delta^\vee$ by a formal integral which would define the invariant linear form once
 the convergence is insured. Finally in \cite{FH}, for $n=1$, the convergence of this linear form is obtained for $\Delta^\sigma|\ |_K^s \times \Delta^\vee|\ |_K^{-s}$ 
and $s$ of real part large enough when $\Delta$ is unitary, the conclusion follows from an analytic continuation argument.\\
We generalize this method here.
The first section is about notations and basic concepts used in the rest of the work.\\
In the second section, we state a theorem of Bernstein (Theorem \ref{Ber}) about rationality of solutions of polynomial systems, and use it as in \cite{CP} or \cite{Ba},
 in order to show, in Proposition \ref{rat}, the holomorphy of integrals of Whittaker functions depending on several complex variables.\\
The third section is devoted to the proof of theorem \ref{dist1}, which asserts that the representation  $\Delta^\sigma|\ |_K^s \times \Delta^\vee|\ |_K^{-s}$ is
 distinguished when $\Delta$ is unitary and $Re(s)$ is in a neighbourhood of $n$.\\
In the fourth section, we extend the result in Theorem \ref{dist2} to every complex number $s$. Our proof relies decisively on a theorem of Youngbin Ok 
(Proposition \ref{bernok} of the present paper), which is a twisted version of a well-known theorem of Bernstein (\cite{Ber}, Theorem A).\\
We end this introduction by recalling a conjecture about classification of distinguished generic representations:

\begin{conj}
Let $m$ be a positive integer, and $\rho$ a generic representation of the group $GL(m,K)$, obtained by normalized parabolic induction of quasi-square-integrable representations $\Delta_1, \dots ,\Delta_t$. It is distinguished if and only if there exists a reordering of the ${\Delta _i}$'s, and an integer $r$ between $1$ and $t/2$, such that we have $\Delta_{i+1}^{\sigma} = \Delta_i^{\vee} $ for $i=1,3,..,2r-1$, and $\Delta_{i}$ is distinguished for $i > 2r$.
\end{conj}

We denote by $\eta$ the nontrivial character of $F^*$ trivial on the norms of $K^*$.
According to Proposition 26 in \cite{F1}, Proposition 12 of \cite{F2}, Theorem 6 of \cite{K}, and Corollary 1.6 \cite{AKT}, our result reduces the proof of the conjecture to show that representations of the form $\Delta_1\times \dots \times \Delta_t$ with $\Delta_{i+1}^{\sigma} = \Delta_i^{\vee} $ for $i=1,3,..,2r-1$ for some $r$ between $1$ and $t/2$, and non isomorphic distinguished or $\eta$-distinguished $\Delta_{i}$'s for $i > 2r$ are not distinguished whenever one of the $\Delta_{i}$'s is $\eta$-distinguished for $i > 2r$. According to \cite{M3}, the preceding conjecture implies the equality of the analytically defined Asai $L$-function and the Galois Asai $L$-function of a generic representation.

\section{Notations}

We denote by $|\ |_K$ and $|\ |_F$ the respective absolute values on $K^*$, by $q_K$ and $q_F$ the respective cardinalities of their residual field, 
by $R_K$ the valuation ring of $K$, and by $P_K$ the maximal ideal of $R_K$. The restriction of  $|\ |_K$ to $F$ is equal to $|\ |_F^2$.\\
More generally, if the context is clear, we denote by $|M|_K$ and $|M|_F$ the positive numbers $|det(M)|_K$ and $|det(M)|_F$ for $M$ a square matrix with determinant in $K$ and 
$F$ respectively. We denote by $G_n$ the algebraic group $GL(n)$.
Hence if $\pi$ is a representation of $G_n(K)$ for some positive $n$, and if $s$ is a complex number, we denote by $\pi|\ |_K^s$ the twist of $\pi$ by the character 
$|det(\ )|_K^s$.\\
We call partition of a positive integer $n$, a family $\bar{n}=(n_1,\dots,n_t)$ of positive integers (for a certain $t$ in $\mathbb{N}-\left\lbrace 0\right\rbrace $), 
such that the sum $n_1+\dots+n_t$ is equal to $n$. To such a partition, we associate a subgroup of $G_n(K)$ denoted by $P_{\bar{n}}(K)$, given by matrices of the form 
$$\left (\begin{array}{cccccccc}
g_1 & \star & \star & \star & \star \\
    & g_2 & \star &  \star & \star \\
    &  & \ddots & \star & \star \\  
    &  &     & g_{t-1} &  \star \\
    &  &     &   & g_t
 \end{array}\right),$$ with $g_i$ in $G_{n_i}(K)$ for $i$ between $1$ and $t$. We call it the standard parabolic subgroup associated with the partition $\bar{n}$. 
We denote by $N_{\bar{n}}(K)$ its unipotent radical subgroup, given by the matrices $$\left (\begin{array}{cccccccc}
I_{n_1} & \star & \star\\
     &  \ddots & \star\\  
    &  & I_{n_t}
 \end{array}\right),$$ and we denote it by $N_n(K)$ when $\bar{n}=(1,\dots,1)$. We denote by $M_{\bar{n}}(K)$ the standard Levi subgroup of matrices 
$\left (\begin{array}{cccccccc}
g_1 &        &   \\  
    & \ddots &    \\
    &        & g_t
 \end{array}\right),$ with $g_i$ in $G_{n_i}(K)$ for $i$ between $1$ and $t$.\\
Finally we denote by $P_n(K)$ the affine subgroup of $GL(n,K)$ given by the matrices $\left (\begin{array}{cccccccc}
g & \star\\  
    & 1
 \end{array}\right)$, with $g$ in $GL({n-1},K)$.\\
Let $X$ be a locally closed space of an $l$-group $G$, and $H$ closed subgroup of $G$, with $H.X \subset X$.
If $V$ is a complex vector space, we denote by $C^{\infty}(X,V)$ the space of smooth functions from $X$ to $V$, and by $C_c^{\infty}(X,V)$ the space of smooth functions with compact support from $X$ to $V$ (if one has $V=\mathbb{C}$, we simply denote it by $C_c^{\infty}(X)$).\\ 
If $\rho$ is a complex representation of $H$ in  $V_{\rho}$, we denote by $C^{\infty}(H \backslash X, \rho, V_{\rho})$ the space of functions $f$ from $X$ to $V_{\rho}$, fixed under the action by right translation of some compact open subgroup $U_f$ of $G$, and which verify $f(hx)=\rho(h)
 f(x)$ for $h \in H$, and $x \in X$ (if $\rho$ is a character, we denote this space by $C^{\infty}(H \backslash X, \rho)$. We denote by $C_c^{\infty}(H \backslash X, \rho, V_{\rho})$ subspace of functions with support compact modulo $H$ of $C^{\infty}(H \backslash X, \rho, V_{\rho})$.\\ 
We denote by $Ind _H^G (\rho)$ the representation by right translation of $G$ in $C^{\infty}(H \backslash G, \rho,
 V_{\rho})$ and by $ind _H^G (\rho)$ the representation by right translation of $G$ in $C_c^{\infty}(H \backslash G,  \rho,
 V_{\rho})$. We denote by ${Ind'} _H^G (\rho)$ the normalized induced representation $Ind _H^G ((\Delta _G /\Delta _H)^{1/2} \rho)$ and by ${ind'} _H^G (\rho)$ the normalized induced representation $ind _H^G ((\Delta _G /\Delta _H)^{1/2} \rho)$.\\ 
Let $n$ be a positive integer, and $\bar{n}=(n_1,\dots,n_t)$ be a partition of $n$, and suppose that we have a representation $(\rho_i, V_i)$ of $GL({n_i},K)$ for each $i$ between $1$ and $t$. Let $\rho$ be the extension to $P_{\bar{n}}(K)$ of the natural representation $\rho_1 \otimes \dots \otimes \rho_t$ of $GL({n_1},K) \times \dots \times GL({n_t},K)$, by taking it trivial on $N_{\bar{n}}(K)$. We denote by $\rho_1 \times \dots \times \rho_t$ the representation ${Ind'} _{P_{\bar{n}}(K)}^{GL(n,K)} (\rho)$.\\

\section{Analytic continuation of Whittaker forms}

If $\rho$ is a generic representation of $G_n(K)$, and $\psi$ is a nontrivial character of $K$, trivial on $F$, then for every $W$ in the Whittaker model 
$W(\rho,\psi)$ of $\rho$, by standard arguments, the following integral
 is convergent for $Re(s)$ large, and defines a rational function in $q_F^{-s}$:

$$I_{(0)}(W,s)=\int_{N_n(F)\backslash P_n(F)} W(p){|det(p)|_F}^{s-1} dp.$$

By standard arguments again, the vector space generated by the
functions $I_{(0)}(W,s)$, for $W$ in $W(\rho,\psi)$, is a fractional ideal $I_{(0)} (\pi)$ of
 $\mathbb{C}[q_F^{-s},q_F^s]$, which has a unique generator which is an Euler factor, independent of $\psi$, that we denote by $L_{F,(0)}^K(\rho,s)$.\\
Similarly, if $\rho'$ is another generic representation of $G_n(K)$, then for every $W$ and $W'$ in the Whittaker models $W(\rho,\psi)$ and  $W(\rho',\psi^{-1})$,
 the following integral
 is convergent for $Re(s)$ large, and defines a rational function in $q_K^{-s}$, which has a Laurent series development in $q_K^{-s}$:

$$I_{(0)}(W,W',s)=\int_{N_n(K)\backslash P_n(K)} W(p)W'(p){|det(p)|_K}^{s-1} dp.$$

The vector space generated by the
functions $I_{(0)}(W,W',s)$, is a fractional ideal of
 $\mathbb{C}[q_K^{-s},q_K^s]$, which has a unique generator which is an Euler factor, independent of $\psi$, that we denote by $L_{(0)}(\rho\times \rho',s)$.\\

According to theorem $9.7$ of \cite{Z}, there is a partition of $n$ and quasi-square-integrable representations $\Delta_1, \dots ,\Delta_t $ associated to it such that 
$\rho$ is isomorphic to $\Delta_1\times \dots \times \Delta_t$. The map $u=(u_1,\dots,u_t)\mapsto q_K^u=(q_K^{u_1},\dots,q_K^{u_t})$ defines an isomorphism of varieties 
between $(\mathcal{D}_K)^t=(\mathbb{C}/\frac{2i \pi}{ln(q_K)\mathbb{Z}})^t$ and $(\mathbb{C}^*)^t$. We also denote by $\mathcal{D}_F$ the variety 
$(\mathbb{C}/\frac{2i \pi}{ln(q_F)\mathbb{Z}})$ which the isomorphism $s\mapsto q_F^{-s}$ identifies to $(\mathbb{C}^*)^t$, and we denote by $\mathcal{D}$ the product 
$(\mathcal{D}_K)^t \times \mathcal{D}_F$.\\
Associate to $u$ and $\rho$ is the representation $\rho_u= \Delta_1 |\ |_K^{u_1}\times \dots \times \Delta_t |\ |_K^{u_t}$. In their classical model, for every representation 
$\rho_u$, the restrictions of the functions of the space of $\rho_u$ to the maximal compact subgroup $GL(n,R_K)$ of $GL(n,K)$ define the same space $\mathcal{F}_{\rho}$, which 
is called the space of flat sections of the series $\rho_u$. To each $f$ in $\mathcal{F}_{\rho}$, corresponds a unique function $f_u$ in $\rho_u$. It is known that for fixed 
$g$ in $GL(n,K)$ and $f$ in $\mathcal{F}_{\rho}$, the function $(u,s)\mapsto |g|_K^s\rho_u(g)f$ belongs to $\mathbb{C}[\mathcal{D}] \otimes_{\mathbb{C}} \mathcal{F}_{\rho}$. 
For every $f$ in $\mathcal{F}_{\rho}$ and $u$ in $(\mathcal{D}_K)^t$, there is a function $W_{f,u}=W_{f_u}$ defined in Section 3.1 of \cite{CP} in the Whittaker model 
$W(\rho_u,\psi)$, such that $W_{f,u}$ describes $W(\rho_u,\psi)$ when $f$ describes $\mathcal{F}_{\rho}$. The space $W^{(0)}$ is defined in \cite{CP} as the complex vector 
space generated by the functions $(g,u)\mapsto W_{f,u}(gg')$ for $g'$ in $GL(n,K)$.\\

We will need a theorem of Bernstein insuring rationality of solutions of polynomial systems. The setting is the following.\\
Let $V$ be a complex vector space of countable dimension. Let $R$ be an index set, and let $\Xi$ be a collection 
$\left\lbrace (x_r, c_r)| r\in R \right\rbrace$ with $x_r \in V$ and $c_r\in \mathbb{C}$. A linear form $\lambda$ in $V^*= Hom_{\mathbb{C}}(V,\mathbb{C})$ is said to be 
a solution of the system $\Xi$ if $\lambda(x_r)=c_r$ for all $r$ in $R$.\\ 
Let $\mathcal{D}$ be an irreducible algebraic variety over $\mathbb{C}$, and suppose that to each $d$, a system 
$\Xi_d =\left\lbrace (x_r(d), c_r(d))| r\in R \right\rbrace$ with the index set $R$ independent of $d$ in $\mathcal{D}$. We say that the family of systems 
$\left\lbrace \Xi_{d}, d \in \mathcal{D}\right\rbrace $ is polynomial if $x_r(d)$ and $c_r(d)$ belong respectively to 
$ \mathbb{C}[\mathcal{D}]\otimes_{\mathbb{C}} V$ and $\mathbb{C}[\mathcal{D}]$. Let $\mathcal{M}=\mathbb{C}(\mathcal{D})$ be the field of fractions of 
$\mathbb{C}[\mathcal{D}]$, we denote by $V_{\mathcal{M}}$ the space $\mathcal{M} \otimes_{\mathbb{C}} V$ and by $V_{\mathcal{M}}^*$ the space 
$Hom_{\mathcal{M}}(V_{\mathcal{M}},\mathcal{M})$.\\
 The following statement is a consequence of Bernstein's theorem, the discussion preceding it, and its corollary in Section 1 of \cite{Ba}.

\begin{thm}{(Bernstein)}\label{Ber}
Suppose that in the above situation, the variety $\mathcal{D}$ is nonsingular and that there exists a non-empty subset $\Omega \subset \mathcal{D}$ open in 
the usual complex topology of $\mathcal{D}$, such that for each $d$ in $\Omega$, the system $\Xi_d$ has a unique solution $\lambda_{d}$. 
Then the system $\Xi=\left\lbrace (x_r(d),c_r(d))|r\in R \right\rbrace$ over the field $\mathcal{M}=\mathbb{C}(\mathcal{D})$ has a unique solution
 $\lambda(d)$ in $V_{\mathcal{M}}^*$, and $\lambda(d)=\lambda_d$ is the unique solution of $\Xi_d$ on $\Omega$. 
\end{thm}

In order to apply this theorem, we first prove the following proposition.

\begin{prop}\label{conv}
Let $\rho$ be a generic representation of $G_n(K)$, there are $t$ affine linear forms $L_i$, for $i$ between $1$ and $t$, with $L_i$ depending on 
the variable $u_i$, such that if the $L_i(u_i)$'s and $s$ have positive real parts, the integral $I_{(0)}(W, s)=\int_{N_n(F)\backslash P_n(F)} W(p)|det(p)|_F^{s-1} dp$ 
is convergent for any $W$ in $W(\rho_u,\psi)$. 
\end{prop}
\begin{proof}

We recall the following claim, which is proved in the lemma of Section 4 of \cite{F1}.

\begin{claim}
Let $\tau$ be a sub-$P_n(K)$-module of $C^{\infty}(N_n(K)\backslash P_n(K), \psi)$, such that for every $k$ between $0$ and $n$, 
the central exponents of the shifted derivatives $\tau^{[k]}$ (see \cite{Ber} 7.2) are positive (i.e. the central characters of all the irreducible sub quotients of $\tau^{[k]}$
 have positive
 real parts), 
then whenever $W$ belongs to $\tau$, the integral $\int_{N_n(F)\backslash P_n(F)} W(p)dp$ is absolutely convergent.
\end{claim}

Applying this to our situation, and denoting by $e_{\rho }$ the maximal element of the set of central exponents of $\rho $ (see Section 7.2 of \cite{Ber}), we deduce that as soon as $u$ is such that $L_i(u)=u_i-e_{\rho } -1$ has positive real part for $i$ between $1$ and $t$, and as soon as $s$ has positive real part, the integral $\int_{N_n(F)\backslash P_n(F)} W(p)|det(p)|_F^{s-1} dp$ converges for all $W$ in $W(\rho_u,\psi)$.
\end{proof}

We now can prove the following:

\begin{prop}\label{rat}
Let $\rho$ be a generic representation of $GL(n,K)$, for every $f$ in $\mathcal{F}_{\rho}$, the function $I_{(0)}(W_{f,u},s)$ belongs to $\mathbb{C}(q_F^{-u},q_F^{-s})$.
\end{prop}
\begin{proof}
In our situation, the underlying vector space is $V=\mathcal{F}_{\rho}$ and is of countable dimension because $\rho$ is admissible. The invariance property satisfied by the functional $I_{(0)}$, for $Re(s)$ large enouigh, 
 is \begin{equation}\label{inv} I_{(0)}(\rho_u(p)W_{f,u},s)=|det(p)|_F^{1-s}I_{(0)}(W_{f,u},s)\end{equation} for $f$ in $\mathcal{F}_{\rho}$, and $p$ in $P_n(F)$.\\
From the proof of Theorem 1 of \cite{K}, it follows that out of the hyperplanes in $(u,s)$ defined by $c_{\rho_u^{(j)}}(t)=|t|_F^{(n-j)(s-1)}$,
 where $\rho_u^{(j)}$ is the representation of $G_{n-j}(F)$ called the $j$-th derivative of $\rho_u$ (see summary before Proposition 2.3 of \cite{AKT}), 
for $j$ from $1$ to $n$, the space of solutions of equation \ref{inv} is of dimension at most one. If we take a basis of 
$(f_{\alpha})_{\alpha\in A}$ of $\mathcal{F}_{\rho}$, the polynomial family over the irreducible complex variety 
$\mathcal{D}=(\mathcal{D}_K)^t\times \mathcal{D}_F$ of systems $\Xi'_d$, for $d=(u,s)\in \mathcal{D}$ expressing the invariance of $I_{(0)}$ is given by:

$$\Xi'_d  = \left\lbrace \begin{array}{lc} ( \rho_u(p)\rho_u(g_i)f_{\alpha}-|det(p)|_F^{1-s} \rho_u(g_i) f_{\alpha},0), \\ \alpha\in A, p\in P_n(F), g_i \in G_n(K)   
\end{array} \right\rbrace $$
 
Now we define $\Omega$ to be the intersection of the three following subsets of $\mathcal{D}$: 
\begin{itemize}
\item the intersection of the complements of the hyperplanes on which uniqueness up to scalar fails,
\item the intersection of the domains $\left\lbrace Re(L_i(u))>0 \right\rbrace $ and $\left\lbrace Re(s)>0 \right\rbrace $, on which $I_{(0)}(W_{f,u},\phi,s)$ is given by an absolutely convergent integral.
\end{itemize}

The functional $I_{(0)}$ is the unique solution up to scalars of the system $\Xi'$, in order to apply Theorem \ref{Ber}, we add for each $d\in \mathcal{D}$ a normalization equation $E_d$ depending polynomially on $d$. This is done as follows.\\
From Proposition 3.4 of \cite{M3}, if $F$ is a positive function in $C_c^{\infty}(N_n(K)\backslash P_n(K), \psi)$, we choose a $W$ in $W_{\rho}^{(0)}$ such that its restriction to $P_n(K)$ is of the form $W(u,p)=F(p) P(q_K^{\pm{u}})$ for some nonzero $P$ in $\mathcal{P}_0$.
We thus have the equality $I_{(0)}(W,u,s)= \int_{N_n(F)\backslash P_n(F)} F(p)|det(p)|_F^{s-1}dp P(q^{\pm{u}})$. Calling $c$ the constant $r\int_{N_n(F)\backslash P_n(F)} F(p)|det(p)|_F^{s-1}dp$, this latter equality becomes $I_{(0)}(W,u,s)=cP(q_K^{\pm{u}})$.\\
 Now as $W$ is in $W_{\rho}^{(0)}$, it can be expressed as a finite linear combination $W(g,u)= \sum_k \rho_u(g_{\alpha}) W_{f_{\alpha},u}(g)$ for appropriate $g_{\alpha}\in GL(n,K)$. Hence our polynomial family of normalization equations (which is actually independent of $s$) can be written $$E_{(u,s)}=\left\lbrace (\sum_{\alpha}  \rho_u(g_{\alpha}) f_{\alpha}, cP(q_K^{\pm{u}}) \right\rbrace.$$ 
We now call $\Xi$ the system given by $\Xi'$ and $E$, it satisfies the hypotheses of Theorem \ref{Ber}, because on the open subset $\Omega$, the functional $I_{(0)}(_,(u,s))$ is well defined and is the unique solution of the system for every $(u,s)$ in $\Omega$. We thus conclude that there is a functional $I'$ which is a solution of $\Xi$ such that $(u,s)\mapsto I'(W_{f,u},s)$ is a rational function of $q_F^{\pm u}$ and $q_F^{\pm s}$ for $f\in \mathcal{F}_{\rho}$. We also know from Theorem \ref{Ber} that $I'(W_{f,u},s)$ is equal to $I_{(0)}(W,u,s)$ on $\Omega$. Hence $I_{(0)}(W,u,s)$ is equal to the rational function $I'(W_{f,u},s)$ when it is defined by a convergent integral for $(u,s)$ in $\Omega$, and we extend it by $I'(W_{f,u},s)$ for general $(u,s)$ (and still denote it by $I_{(0)}(W,u,s)$).\end{proof}

We now recall the following theorem of Youngbin Ok:

\begin{prop}{(\cite{Ok}, Theorem 3.1.2 or Proposition 1.1 of \cite{M2})} \label{bernok}
Let $\rho$ be an irreducible distinguished representation of $G_n(K)$, if $L$ is a $P_n(F)$-invariant linear form on the space of $\rho$, then it is actually 
$G_n(F)$-invariant.
\end{prop}

We also recall the proposition 2.3 of \cite{M2}. 

\begin{prop}\label{invform}
Let $\rho$ be a generic representation of $G_n(K)$, for any $s\in \mathbb{C}$, the functional ${\Lambda}_{\rho,s}: W \mapsto
  I_{(0)}(W,s)/L_{(0)}(\rho,s)$ defines a nonzero linear form on
 $W(\rho,\psi)$ which transforms
 by $| \  |_F^{1-s}$ under the affine subgroup $P_n(F)$.\\
For fixed $W$ in $W(\rho,\psi)$, the function $s \mapsto {\Lambda}_{\rho,s}(W)$ is a polynomial of
 $q_F^{-s}$.

\end{prop}

\section{Distinction of representations $\pi^{\sigma}|\ |_K^s\times \pi^{\vee}|\ |_K^{-s}$ for  unitary $\pi$, Re(s) near n}

We denote by $G$ the group $GL(2n,K)$, by $H$ its subgroup $GL(2n,F)$, by $G'$ the group $GL(n,K)$ and by $M$ the group $M_n(K)$. We denote by $P$ the
group $P_{(n,n)}(K)$, and by $N$ the group $N_{(n,n)}(K)$.\\ 
We denote by $\bar{H}$ subgroup of $G$ given by matrices of the form $\left( \begin{array}{lr}
 A & B \\ B^{\sigma} & A^{\sigma}
\end{array}\right)$, and by $\bar{T}$ the subgroup of $\bar{H}$ of matrices $\left( \begin{array}{lr}
 A & 0 \\ 0 & A^{\sigma}      
 \end{array}\right) $, with $A$ in $G'$.\\
We let $\delta$ be an element of $K-F$ whose square belongs to $F$, and let $U$ be the matrix $\left( \begin{array}{lr}
 I_n & -\delta I_n \\ I_n & \delta I_n     
\end{array}\right)$ of $G$, and $W$ the matrix $\left(\begin{array}{lr}
  & I_n \\                                                              
I_n     &         
\end{array} \right)$. One has $U^{\sigma}U^{-1}=W$ and the group $H$ is equal to $U^{-1} \bar{H} U$.

\begin{LM}\label{openclass}
The double class $PUH$ is opened in $G$.
\end{LM}
\begin{proof}
Call $S$ the space of matrices $g$ in $G$ verifying $g^{\sigma}=g^{-1}$, which is, from Proposition 3. of chapter 10 of \cite{S}, homeomorphic to the quotient space $G / H$ by the map $Q: g\mapsto g^{\sigma}g^{-1}$. As the map $Q$ sends $U$ on $W$, the double class $PUH$ corresponds to the open subset of matrices $\left(\begin{array}{cc} A & B\\ C & D \end{array} \right)$ in $S$ such that $det(C)\neq 0$, the conclusion follows.

\end{proof}

We prove the following integration formula.

\begin{LM}\label{intform}
There is a right invariant measure $d\dot{h}$ on the quotient space $\bar{T}\backslash \bar{H}$, and a Haar measure $dB$ on $M$, such that for any measurable positive function $\phi$ on the quotient space $\bar{T}\backslash \bar{H}$, then the integrals $$\int_{\bar{T}\backslash \bar{H}} \phi(\dot{h})d\dot{h}$$ and  $$\int_{M} \phi\left( \begin{array}{lr}
 I_n & B \\                                                                                                                            B^{\sigma} & I_n        
\end{array}\right) \frac{dB}{|I_n- BB^{\sigma}|_K^n}$$ are equal.

\end{LM}
\begin{proof} It suffices to show this equality when $\phi$ is positive, continuous with compact support in $\bar{T}\backslash \bar{H}$. We fix Haar measures $dt$ on $\bar{T}$ and $dg$ on $\bar{H}$, such that $d\dot{h}dt=dg$. It is known that there exists some positive function $\tilde{\phi}$ with compact support in $\bar{H}$, such that $\phi=\tilde{\phi}^{\bar{T}}$, which means that for any $\dot{h}$ in $\bar{H}$, one has $\phi(\dot{h})=\int_{\bar{T}}\tilde{\phi}(tg) dt$. One then has the relation $$\int_{\bar{T}\backslash \bar{H}} \phi(\dot{h})d\dot{h}= \int_{\bar{H}} \tilde{\phi}(g)dg.$$
Now as $\bar{H}$ is conjugate to $H$, there are Haar measures $dA$ and $dB$ on $M$ such that $dt$ is equal to $d^*A=\frac{dA}{|A|_K^n}$, and the Haar measure on $\bar{H}$ is described by the relation $$d\left( \begin{array}{lr}
 A & B \\                                                                                                             B^{\sigma} & A^{\sigma}        
\end{array}\right)= \frac{dA dB}{\left|\left( \begin{array}{lr}
 A & B \\                                                                                                             B^{\sigma} & A^{\sigma}        
\end{array}\right)\right| _F^{2n}}=\frac{dA dB}{\left|\left( \begin{array}{lr}
 A & B \\                                                                                                             B^{\sigma} & A^{\sigma}        
\end{array}\right)\right| _K^n}.$$  

Hence we have $$\begin{array}{ll} \int_{\bar{T}\backslash \bar{H}} \phi(\dot{h})d\dot{h}& = \int_{M\times M} \tilde{\phi}\left( \begin{array}{lr}
 A & B \\                                                                                                             B^{\sigma} & A^{\sigma}        
\end{array}\right)\frac{dA dB}{\left|\left( \begin{array}{lr}
 A & B \\                                                                                                             B^{\sigma} & A^{\sigma}        
\end{array}\right)\right| _K^n}\\
                              
                                 &= \int_{M\times M} \tilde{\phi} \left[ \left(\begin{array}{lr}
 A &  \\                                                                                                                & A^{\sigma}         
\end{array} \right) \left( \begin{array}{lr}
 I_n & A^{-1}B \\                                                                                                    (A^{-1}B)^{\sigma} & I_n        
\end{array}\right)\right] \frac{dA dB}{|A|_K^{2n} |I_n - A^{-1}B (A^{-1}B)^{\sigma}|_K^n} \end{array}$$ as the complement of $G'$ is a set of measure zero of $M$ (we recall that if $M$ is in $G'$, one has $det \left( \begin{array}{lr}
 I & M \\                                                                                                             M^{\sigma} & I        
\end{array}\right)\!\! = \!det\left(\left(\!\!\! \begin{array}{lr}
 I & M \\                                                                                                             M^{\sigma} & I        
\end{array}\!\!\!\right)\left(\!\!\! \begin{array}{lr}
 I &  \\                                                                                                     - M^\sigma & I        
\end{array}\!\!\!\right)\right)\!\! =  \!det \left(\!\!\! \begin{array}{lr}
 I-MM^\sigma & M \\                                                                                                       & I        
\end{array}\!\!\!\right)\!\! =  \!det(I-MM^\sigma)$).

This becomes after the change of variable $B:=A^{-1}B$ equal to $$\int_{M\times M} \tilde{\phi} \left[ \left(\begin{array}{ll}
 A &  \\                                                                                                                & A^{\sigma}         
\end{array} \right) \left( \begin{array}{lr}
 I_n &  B \\                                                                                                    B^{\sigma} & I_n        
\end{array}\right)\right] \frac{dA}{|A|_K^n} \frac{dB}{|I_n - BB^{\sigma}|_K^n}$$ which is itself equal to $$\int_{G'\times M} \tilde{\phi} \left[ \left(\begin{array}{lr}
 A &  \\                                                                                                                & A^{\sigma}         
\end{array} \right) \left( \begin{array}{lr}
 I_n &  B \\                                                                                                    B^{\sigma} & I_n        
\end{array}\right)\right] d^*A \frac{dB}{|I_n - BB^{\sigma}|_K^n}.$$
The conclusion follows from the fact that $\tilde{\phi}^{\bar{T}}$ is equal to $\phi$.\end{proof}

\begin{thm}\label{dist1}
Let $n$ be a positive integer, and let $\pi$ be a generic unitary representation of $G'$. Then the representation $\pi^{\sigma}|\ |_K^{s} \times \pi^{\vee}|\ |_K^{-s}$ is a 
distinguished representation of 
$G$ for $s$ of real part in a neighbourhood of $n$.
\end{thm}

\begin{proof}

We denote by $\Pi_s$ the representation $\pi^{\sigma}|\ |_K^{s} \times \pi^{\vee}|\ |_K^{-s}$ of $G_{2n}(K)$. Let $V$ be the space of the representation 
$\pi$ (and $\pi^\sigma$), and $V^\vee$ be the space of its smooth contragredient $\pi^\vee$.

We first start with the following lemma: 

\begin{LM}\label{bound}
Any coefficient of the representation $\Pi_u$ is bounded for $Re(u)$ near zero.
\end{LM}

\begin{proof}[Proof of Lemma \ref{bound}]
From \cite{Ber}, as $\Pi_0$ (resp. $\Pi_0^\vee$) is unitary, we know that all its shifted derivatives have positive central exponents. Actually,
 this latter property remains true for
 $\Pi_u$ (resp. $\Pi_u^\vee$) for $Re(u)$ in a neighbourhood of zero.\\
Realizing $\Pi_u$ in its Whittaker $ W(\Pi_u,\psi)$, it is a consequence of \cite{Ber}, Theorem B (since $\Pi_u$ is always irreducible for $Re(u)$ near zero), and of 
Proposition \ref{l2} of the appendix, 
that any coefficient of $\Pi_u$ is of the form $g\mapsto \int_{N_n(K)\backslash P_n(K)} \Pi_u(g) W(p) W'(p) dp$, for $W$ in $ W(\Pi_u,\psi)$ and $W'$ in 
$ W(\Pi_u^\vee,\psi^{-1})$.  But then, from Cauchy-Schwartz inequality, we see that any coefficient of $\Pi_u$ is bounded on $P_n(K)$. Now because the central character of 
$\Pi_u$ is unitary, this implies that any coefficient of $\Pi_u$ is bounded on the maximal parabolic subgroup $P_{n-1,1}(K)$. This latter fact, combined with the 
Iwasawa decomposition in $G_n(K)$, and the smoothness of the coefficients, implies that the coefficients of $\Pi_u$ are actually bounded on $G_n(K)$.

\end{proof}

 We denote by $L$ the linear form on $V\otimes V^\vee$ who sends the elementary tensor $v\otimes v^\vee$ to $v^\vee(v)$, it is clearly invariant under the group 
$\pi^\sigma \otimes \pi^\vee (\bar{T})$.\\

Step 1.\\
 We denote by $\rho_s$ the representation $P$, which is the extension of $\pi^{\sigma}|\ |_K^s \otimes \pi^{\vee}|\ |_K^{-s}$ to $P$ by 
the trivial representation of $N_{(n,n)}(K)$.
Here for every $s$, the group $P$ acts through the representation $\rho_s$ on $V\otimes V^\vee$.\\
As a function $f_s$ in the space $C^{\infty}_c(P\backslash G,\Delta_P^{-1/2} \rho_s)$ of $\pi^{\sigma}|\ |_K^s \times \pi^{\vee}|\ |_K^{-s}$, verifies relation 
$$ f_s \left[ \left( \begin{array}{lr}
 A & \star \\                                                                                                                            0 & B        
\end{array}\right)g\right] = \frac{|det(A)|_K^{n/2+s}}{|det(B)|_K^{n/2+s}} \pi^{\sigma}(A) \otimes \pi^{\vee}(B) f(g),$$
we deduce that the restriction to $\bar{H}$ of the function $L_{f_s}: g \mapsto L(f_s(g))$ belongs to the space $C^{\infty} (\bar{T}\backslash \bar{H})$, but its support modulo $\bar{T}$ is
 generally not compact, we will show later that the space of functions obtained this way contains $C^{\infty}_c (\bar{T}\backslash \bar{H})$ as a proper subspace. We must show 
that for $s=n+u$ with $u$ near zero, the integral $\int_{\bar{T}\backslash \bar{H}}|L_{f_s}(\dot{h})|d\dot{h}$ converges.\\
Denoting by $\eta_s$ the function on $G'$ defined by $\eta_s[\left(\begin{array}{lr} A_1 & X\\  & A_2\end{array}\right) k]=(\frac{ |A_1|_K }{|A_2|_K})^{s}$, 
For any complex numbers $t$ and $u$, the multiplication map 
$f_u\mapsto f_{t+u}= \eta_t f_u$ is a vector space isomorphism between $C^{\infty}_c(P\backslash G,\Delta_P^{-1/2} \rho_u)$ and 
$C^{\infty}_c(P\backslash G,\Delta_P^{-1/2} \rho_{t+u})$.\\

In the following, all equalities will be formal, we will show that they have sense for $t=n$ and $u$ near zero at the end of this step, by proving the absolute
 convergence of the considered integrals.
According to lemma \ref{intform}, the integral $\int_{\bar{T}\backslash \bar{H}}|L_{f_{t+u}}(\dot{h})|d\dot{h}$ is equal to $$\int_{M} |L_{f_{t+u}}|\left( \begin{array}{lr}
 I_n & B \\                                                                                                          B^{\sigma} & I_n        
\end{array}\right) \frac{dB}{|I_n- BB^{\sigma}|_K^n}= \int_{M} \eta_{Re(t)}\left( \begin{array}{lr}
 I_n & B \\                                                                                                          B^{\sigma} & I_n        
\end{array}\right) |L_{f_u}|\left( \begin{array}{lr}
 I_n & B \\                                                                                                          B^{\sigma} & I_n        
\end{array}\right) \frac{dB}{|I_n- BB^{\sigma}|_K^n}.$$ 

Now we suppose that $Re(u)$ is near zero. We remind that the quantity $|L_{f_u}|\left( \begin{array}{lr}
 I_n & B \\                                                                                                          B^{\sigma} & I_n        
\end{array}\right)$ is defined for $B$ such that $det\left( \begin{array}{lr}
 I_n & B \\                                                                                                          B^{\sigma} & I_n        
\end{array}\right)\neq 0$, we claim that it is actually bounded by some positive real number $M$.\\
 Indeed the linear for $v^\vee: f_u\mapsto L(f_u(I_{2n}))$ belongs to the smooth dual of $\Pi_u$, and the
 coefficient $|L_{f_u}|(g)$ which equals $|<v^\vee,\Pi_u(g)f_u>|$, is bounded by Lemma \ref{bound}\\

As before, we can suppose that $B$ belongs to $G'$, hence the following decomposition holds
$$\left (\begin{array}{lr}
 I_n & B \\                                                                                                           B^{\sigma}     & I_n         
\end{array} \right)= \left(\begin{array}{lr}
 (-I_n + BB^{\sigma})B^{-\sigma} & I_n \\                                                                                                               & B^{\sigma}          
\end{array} \right)\left(\begin{array}{lr}
  & I_n \\                                                                                                           I_n     &         
\end{array} \right)\left(\begin{array}{lr}
 I_n & B^{-\sigma} \\                                                                                                                & I_n         
\end{array} \right).$$

Denoting by $\tilde{\eta}_s$ the function $g\mapsto \eta_s( wg)$, we only need to look at the convergence of the integral:

$$\begin{array}{ll}
 \int_{M} \eta_{Re(t)} \left( \begin{array}{lr}
 I_n & B \\                                                                                                          B^{\sigma} & I_n        
\end{array}\right) \frac{dB}{|I_n- BB^{\sigma}|_K^n}  &= \int_{M} (\frac{|BB^{\sigma} - I_n|_K}{|B|_K^2}) ^{Re(t)} \tilde{\eta}_{Re(t)} \left(\begin{array}{lr}
 I_n & B^{-\sigma} \\                                                                                                                & I_n         
\end{array} \right)\frac{dB}{| I_n- BB^{\sigma}|_K^n} \\
                   &=\int_{G'} (\frac{| I_n- BB^{\sigma}|_K}{|B|_K^2}) ^{Re(t)}\tilde{\eta}_{Re(t)}\left(\begin{array}{lr}
 I_n & B^{-\sigma} \\                                                                                                                & I_n         
\end{array} \right) \frac{|B|_K^n d^*B}{|I_n- BB^{\sigma}|_K^n}\\ 
                   &=\int_{G'} \frac{(|I_n-  C^{-\sigma} C^{-1}|_K {|C|_K^2})^{Re(t)}}{|I_n-  C^{-\sigma} C^{-1}|_K^n}  \tilde{\eta}_{Re(t)}\left(\begin{array}{lr}
 I_n & C \\                                                                                                                & I_n         
\end{array} \right) \frac{d^* C}{|C|_K^n} \\
                   &=\int_{M} |C C^{\sigma} - I_n|_K^{Re(t) -n} \tilde{\eta}_{Re(t)} \left(\begin{array}{lr}
 I_n & C \\                                                                                                                & I_n         
\end{array} \right) d C
\end{array}$$

We recognize here the function $\tilde{\eta}$ of 4. (3) of \cite{JPS1} (p.411). The following lemma and its demonstration was communicated to me by Jacquet.

\begin{LM}{(Jacquet)} Let $\Phi_0$ be the characteristic function of $M_{n}(R_K)$, then from the Godement-Jacquet
 theory of Zeta functions of simple algebras, the integral $\int_{G'}\phi_0 (H) |H|_K^u d^*H$ is convergent for $Re(t)\geq n-1$, and is equal to $1/P(q_K^{-t})$ 
for a nonzero polynomial $P$. Then, for $Re(t)\geq (n-1)/2$, and $g$ in $G'$, denoting by $\phi$ the characteristic function of $M_{n,2n}(R_K)$
(matrices with $n$ rows and $2n$ columns) one has $$\tilde{\eta}_t(g)=P(q_K^{-2s})|g|_K^t \int_{G'} \Phi[(H,0)g]|H|_K^{2t}d^*H.$$
 \end{LM}
\begin{proof}[Proof of the lemma]
It is a consequence of the decomposition $G'=N_{(n,n)}^- (K) M_{(n,n)}(K) G_{2n}(R_K)$ (with $N_{(n,n)}^- (K)$ the opposite of $N_{(n,n)} (K)$),
and of the fact that functions on both sides verify the relation
 $f[\left(\begin{array}{lr} A_1 &  \\ X & A_2\end{array}\right) g]=\frac{|A_2|_K^t}{|A_1|_K^t} f(g)$, and are both equal to $1$ on $G_{2n}(R_K)$ 
(if $d^*H$ is normalized so that the maximal compact subgroup $G_{2n}(R_K)$ has measure $1$).
\end{proof}

Finally, we suppose moreover that $Re(t)=n$, hence we need to check the convergence of $$\int_{M} \tilde{\eta}_n\left(\begin{array}{lr}
 I_n & C \\  & I_n         
\end{array} \right) d C= \int_{M} (P(q_K^{-2n}) \int_{G'} \Phi[(H,HC)]|H|_K^{2n}d^*H)dC$$

As the functions in the integrals are positive, by Fubini's theorem, this latter is equal to:

$$  P(q_K^{-2s})  \int_{G'} \int_{M} \Phi[(H,HC)]dC |H|_K^{2n}d^*H $$
$$ = P(q_K^{-2s})  \int_{G'} \int_{M} \Phi[(H,C)]dC |H|_K^{n}d^*H $$
$$= P(q_K^{-2s})  \int_{M} \int_{M} \Phi[(H,C)]dC $$ which clearly converges.

Hence we proved that the $H$-invariant linear form $f_{s} \mapsto \int_{\bar{T}\backslash \bar{H}} L_{f_{s}}(\dot{h})d\dot{h}$ on $\Pi_s$ was well defined for $Re(s)$ in a 
neighbourhood of $n$.

Step 2.\\
Suppose that the complex number $s$ has real part greater than $n$. We are going to show that the linear form 
$\Lambda: f_s\mapsto \int_{\bar{T}\backslash \bar{H}}L_{f_s}(\dot{h})d\dot{h}$ is nonzero. More precisely we are going to show that the space of functions 
$L(f)$ on $\bar{T}\backslash \bar{H}$ for $f$ in $C_c^{\infty}(P\backslash G, {\Delta}_{P}^{-1/2} \rho_s)$, contain $C^{\infty}_c (\bar{T}\backslash \bar{H})$.\\
According to Lemma \ref{openclass} , the double class $PUH$ is opened in $G$, hence the extension by zero outside $PUH$ gives an injection of the space
 $C_c^{\infty}(P \backslash PUH, {\Delta}_{P}^{-1/2}\rho_s)$ into the space $C_c^{\infty}(P\backslash G, {\Delta}_{P}^{-1/2} \rho_s)$.\\
But the right translation by $U$, which is a vector space automorphism of $C_c^{\infty}(P\backslash G, {\Delta}_{P}^{-1/2} \rho_s)$, sends 
$C_c^{\infty}(P \backslash PUH, {\Delta}_{P}^{-1/2}\rho_s)$ onto $C_c^{\infty}(P \backslash P\bar{H}, {\Delta}_{P}^{-1/2}\rho_s)$, hence
 $C_c^{\infty}(P \backslash P\bar{H}, {\Delta}_{P}^{-1/2}\rho_s)$ is a subspace of $C_c^{\infty}(P\backslash G, {\Delta}_{P}^{-1/2} \rho_s)$.\\ 
Now restriction to $\bar{H}$ defines an isomorphism between
 $C_c^{\infty}(P \backslash P\bar{H}, {\Delta}_{P}^{-1/2}\rho_s)$ and $C^{\infty}_c (\bar{T}\backslash \bar{H}, \rho_s)$ because ${\Delta}_{P}$ has trivial restriction 
to the group $\bar{T}$. 
But then the map $f\mapsto L(f)$ defines a morphism of $\bar{H}$-modules from $C^{\infty}_c (\bar{T}\backslash \bar{H}, \rho_s)$ to 
$C^{\infty}_c (\bar{T}\backslash \bar{H})$, which is surjective because of the commutativity of the following diagram,

$$\begin{array}{ccc}
C^{\infty}_c (\bar{H})\otimes V_{\rho_s} & \overset{Id\otimes L}{\longrightarrow} & C^{\infty}_c (\bar{H})\\
\downarrow  &  & \downarrow \\
C^{\infty}_c (\bar{T}\backslash \bar{H}, \rho_s) & \longrightarrow & C^{\infty}_c (\bar{T}\backslash \bar{H})
\end{array},$$
where the vertical arrows defined in Lemma 2.9 of \cite{M1} and the upper arrow are surjective.\\
We thus proved that space of restrictions to $\bar{H}$ of functions of $L(f)$, for $f$ in $C_c^{\infty}(P\backslash G, {\Delta}_{P}^{-1/2} \rho_s)$, 
contain $C^{\infty}_c (\bar{T}\backslash \bar{H})$, hence $\Lambda$ is nonzero and the representation $\pi^{\sigma}|\ |_K^s \times \pi^{\vee}|\ |_K^{-s}$
is distinguished for $Re(s)$ near $n$.\\

\end{proof}

\section{Distinction of $\Delta^{\sigma}\times \Delta^{\vee}$ for quasi-square-integrable $\Delta$}

Now we are going to restrain ourself to the case of $\pi$ a discrete series representation.\\
We recall  if $\rho$ is a supercuspidal representation of $G_r(K)$ for a positive integer $r$. The representation $\rho \times \rho|\ |_F \times \dots \times \rho|\ |_F^{l-1}$
 of $G_{rl}(K)$ is reducible, with a unique irreducible quotient that we denote by $[\rho|\ |_K^{l-1},\rho|\ |_K^{l-2},\dots, \rho]$. A representation $\Delta$ of the group 
$G_n(K)$ is quasi-square-integrable if and only if there is $r\in \left\lbrace 1,\dots, n \right\rbrace$ and $l\in \left\lbrace 1,\dots, n \right\rbrace$ with $lr=n$, and 
$\rho$ a supercuspidal representation of $G_r(K)$ such that the representation $\Delta$ is equal to $[\rho|\ |_K^{l-1},\rho|\ |_K^{l-2},\dots, \rho]$, the representation 
$\rho$ is unique.\\
Let $\Delta_1$ and $\Delta_2$ be two quasi-square-integrable representations  of $G_{l_1 r}(K)$ and $G_{l_2 r'}(K)$, of the form 
$[\rho_1|\ |_K^{l_1-1},\rho_1|\ |_K^{l_1-2},\dots, \rho_1]$ with $\rho_1$ a supercuspidal representation of $G_{r}(K)$, and\\ 
$[\rho_2|\ |_K^{l_2-1},\rho_1|\ |_K^{l_2-2},\dots, \rho_2]$  with $\rho_2$ a supercuspidal representation of $G_{r'}(K)$ respectively, then if 
$\rho_1=\rho_2|\ |_K^{l_2}$, we denote by $[\Delta_1,\Delta_2]$ the quasi-square integrable representation $[\rho_1|\ |_K^{l_1-1},\dots,\rho_2]$ of $G_{(l_1+l_2) r}(K)$.
Two quasi-square-integrable representations $\Delta=[\rho|\ |_K^{l-1},\rho|\ |_K^{l-2},\dots, \rho]$ and $\Delta'=[\rho'|\ |_K^{l'-1},\rho'|\ |_K^{l'-2},\dots, \rho']$ of 
$G_n(K)$ and $G_{n'}(K)$ are said to be linked if  $\rho'=\rho|\ |_K^{k'}$ with $k'$ between $1$ and $l$, and $l'>l$, or if $\rho=\rho'|\ |_K^k$, with $k$ between $1$ 
and $l'$, and $l>l'$. It is known that the representation $\Delta \times \Delta'$ always has a nonzero Whittaker functional on its space, and is irreducible if and only if 
$\Delta$ and $\Delta'$ are unlinked.\\
 
We will need the following theorem.

\begin{thm}\label{cndist}
Let $n_1$ and $n_2$ be two positive integers, and $\Delta_1$ and $\Delta_2$ be two unlinked quasi-square integrable representations of $G_{n_1}(K)$ and $G_{n_2}(K)$ 
respectively. If the representation $\Delta_1\times \Delta_2$ of $G_{n_1+n_2}(K)$ is distinguished, then either both  $\Delta_1$ and $\Delta_2$ are distinguished, either
 $\Delta_2^\vee$ is isomorphic to $\Delta_1^\sigma$.
\end{thm}

\begin{proof} In the proof of this theorem, we will denote by $G$ the group $G_{n_1+n_2}(K)$ (not the group $G_{2n}(K)$ anymore), by $H$ the group $G_{n_1+n_2}(F)$, and by $P$ the group $P_{(n_1,n_2)}(K)$.\\
As the representation  $\Delta_1\times \Delta_2$ is isomorphic to $\Delta_2\times \Delta_1$, we suppose $n_1\leq n_2$.
From Lemma 4 of \cite{F4}, the $H$-module $\pi$ has a factor series with factors isomorphic to the representations $ind_{u^{-1}Pu\cap H}^H((\delta_P^{1/2}\Delta_1\otimes\Delta_2 )^u)$ (with $(\delta_P^{1/2}\Delta_1\otimes\Delta_2 )^u (x)=\delta_P^{1/2}\Delta_1\otimes\Delta_2 (uxu^{-1})$) when $u$ describes a set of representatives of $P\backslash G/H$. 
Hence we first describe such a set. 

\begin{LM}\label{repr}
The matrices $u_k=\begin{bmatrix}
I_{n_1-k} &  &  &  \\
          & I_k & -\delta I_k &  \\
          & I_k & \delta I_k &  \\
    &  &  & I_{n_2-k} \end{bmatrix}$, give a set of representatives $R(P\backslash G/H)$ of the double classes $P\backslash G/H$ when $k$ describes the set $\{0,\dots,n_1\}$ (we set $u_0=I_{n_1+n_2}$). 
 \end{LM}
\begin{proof}[Proof of Lemma \ref{repr}]
Set $n=n_1+n_2$, the quotient set $H\backslash G/P$ identifies with the set of orbits of $H$ for its action on the variety of $K$-vectors spaces of dimension $n_1$ in $K^{n}$. We claim that two vector subspaces $V$ and $V'$ of dimension $n_1$ of $K^n$ are in the same $H$-orbit if and only if $dim(V\cap V^\sigma)$ equals $dim(V'\cap {V'}^\sigma)$. This condition is clearly necessary. If it is verified, we choose $S$ a supplementary space of $V\cap V^\sigma$ in $V$ and we choose $S'$ a supplementary space of $V'\cap {V'}^\sigma$ in $V'$, $S$ and $S'$ have same dimension. We also choose $Q$ a supplementary space of $V+V^\sigma$ in $K^n$ defined over $F$ (i.e. stable under $\sigma$, or equivalently having a basis in the space $F^n$ of fixed points of $K^n$ under $\sigma$), and $Q'$ a supplementary space of $V'+{V'}^\sigma$ in $K^n$ defined over $F$, and $Q$ and $Q'$ have the same dimension. Hence we can decompose $K^n$ in the two following ways: $K^n= (V\cap V^\sigma) \oplus (S\oplus S^\sigma) \oplus Q$ and $K^n= (V'\cap {V'}^\sigma) \oplus (S'\oplus {S'}^\sigma) \oplus Q'$. Let $u_1$ be an isomorphism between $V\cap V^\sigma$ and $V'\cap {V'}^\sigma$ defined over $F$ (i.e. $u(v_1^\sigma)=u(v_1)^\sigma$ for $v_1$ in $V\cap V^\sigma$), $u_2$ an isomorphism between $S$ and $S'$ (to which we associate an isomorphism $u_3$ between $S^\sigma$ and ${S'}^\sigma$ defined by $u_3(v)=(u_2(v^\sigma))^\sigma$ for $v$ in $S^\sigma$), and $u_4$ an isomorphism between $Q$ and $Q'$ defined over $F$. Then the isomorphism $h$ defined by $v_1+v_2+v_3+v_4\mapsto u_1(v_1)+u_2(v_2)+u_3(v_3)+u_4(v_4)$ is defined over $F$, and sends $V=S\oplus V\cap V^\sigma$ to $V'=S'\oplus V'\cap {V'}^\sigma$, hence $V$ and $V'$ are in the same $H$-orbit.\\
 If $(e_1,\dots,e_n)$ is the canonical basis of $K^n$, we denote by $V_{n_1}$ the space $Vect(e_1,\dots,e_{n_1})$. Let $k$ be an integer between $0$ and $n_1$, the image $V_k$ of $V_{n_1}$ by the morphism whose matrix in the canonical basis of $K^n$ is $\begin{bmatrix}
I_{n_1-k} &  &  &  \\
          & 1/2 I_k & 1/2 I_k &  \\
          & -1/(2\delta) I_k & 1/(2\delta) I_k &  \\
    &  &  & I_{n_2-k} \end{bmatrix}$ verifies $dim(V_{n_1}\cap V_{n_1}^\sigma)=n_1-k$. Hence the matrices $\begin{bmatrix}
I_{n_1-k} &  &  &  \\
          & 1/2 I_k & 1/2 I_k &  \\
          & -1/(2\delta) I_k & 1/(2\delta) I_k &  \\
    &  &  & I_{n_2-k} \end{bmatrix}$ for $k$ between $0$ and $n_1$ give a set of representatives of the quotient set $H\backslash G/P$, which implies that their inverses $\begin{bmatrix}
I_{n_1-k} &  &  &  \\
          & I_k & -\delta I_k &  \\
          & I_k & \delta I_k &  \\
    &  &  & I_{n_2-k} \end{bmatrix}$ give a set of representatives of $P\backslash G/H$.

\end{proof}

We will also need to understand the structure of the group $P\cap uHu^{-1}$ for $u$ in $R(P\backslash G/H)$. 

\begin{LM}\label{structure}
Let $k$ be an integer between $0$ and $n_2$, we deduce the group $P\cap u_k Hu_k^{-1}$ is the group of matrices of the form $\begin{bmatrix}
H_1 &  X  & X^\sigma & M \\ 
          & A &   & Y \\
          &    & A^\sigma & Y^\sigma \\
          &     &   &           H_2 \end{bmatrix}$
 for $H_1$ in $G_{n_1-k}(F)$, $H_2$ in $G_{n_2-k}(F)$, $A$ in $G_k(K)$, $X$ in $M_{n_1 -k, k}(K)$, $Y$ in $M_{k, n_2 -k}(K)$, and $M$ in $M_{n_1-k,n_2-k}(F)$. It is the
 semi-direct product of the subgroup $M_k(F)$ of matrices of the preceding form with $X$, $Y$, and $M$ equal to zero, and of the subgroup $N_k$ of matrices of the preceding 
form with $H_1=I_{n_1-k}$, $H_2=I_{n_2-k}$, and $A=I_k$. Moreover denoting by $P_k$ the parabolic subgroup of $M_{(n_1,n_2)}(K)$ associated with the sub partition
 $(n_1-k,k,k,n_2-k)$ of $(n_1,n_2)$, the following relation of modulus characters is verified: ${\delta^2_{P\cap u_k Hu_k^{-1}|M_k(F)}}= (\delta_{P_k}\delta_{P})_{|M_k(F)}$. 
\end{LM}

\begin{proof}[Proof of Lemma \ref{structure}]
One verifies that the algebra $u_k M_n(K)u_k^{-1}$ consists of matrices having the block decomposition corresponding to the partition $(n_ 1-k,k,k, n_2-k)$ of the form 
$\begin{bmatrix}
M_1       &  X  & X^\sigma  & M_2 \\
Y         &  A  &  B^\sigma &    Y' \\
Y^\sigma  &  B   & A^\sigma & {Y'}^\sigma \\
M_3          &  X' &  {X'}^\sigma &  M_4 \end{bmatrix}$, the first part of the proposition follows. For the second part, if the matrix $T=\begin{bmatrix}
H_1 &   &  &  \\
          & A &   &  \\
          &    & A^\sigma &  \\
          &     &   &           H_2 \end{bmatrix}$
belongs to $M_k(F)$, the complex number $ \delta_{P\cap u_k Hu_k^{-1}}(T)$ is equal to the modulus of the automorphism $int_T$ of $N_k$, hence is equal to 
$$|H_1|_F^{2k}|A|_K^{k-n_1}|H_1|_F^{n_2-k}|H_2|_F^{k-n_1}|A|_K^{n_2-k}|H_2|_F^{-2k}=|H_1|_F^{n_2+k}|A|_K^{n_2-n_1}|H_2|_F^{-k-n_1}.$$
 In the same way, the complex number $ \delta_{P_k}(T)$ equals $$|H_1|_K^{k}|A|_K^{k-n_1}|A|_K^{n_2-k}|H_2|_F^{-k}= |H_1|_F^{2k}|A|_K^{n_2-n_1}|H_2|_F^{-2k},$$ and 
$ \delta_{P_k}(T)$ equals $(|H_1|_K|A|_K)^{n_2})(|H_2|_K|A|_K)^{-n_1})=|H_1|_F^{2n_2}|A|_K^{n_2-n_1}|H_2|_F^{-2n_1}$.\\
 The wanted relation between modulus characters follows.
\end{proof}

A helpful corollary is the following.

\begin{cor}\label{util}
Let $P_k$ be the standard parabolic subgroup of $M_{(n_1,n_2)}(K)$ associated with the sub partition $(n_1-k,k,k,n_2-k)$ of $(n_1,n_2)$, $U_k$ its unipotent radical, and 
$N_k$ the intersection of the unipotent radical of the standard parabolic subgroup of $G$ associated with the partition $(n_1-k,k,k,n_2-k)$ and $uHu^{-1}$. Then one has 
$U_k\subset N_k N$.
\end{cor}
\begin{proof}[Proof of Corollary \ref{util}]
It suffices to prove that matrices of the form $\begin{bmatrix}
I_{n_1-k} &  X &  &  \\
          & I_k &   &  \\
          &    & I_k &  \\
          &     &   & I_{n_2-k} \end{bmatrix}$ and $\begin{bmatrix}
I_{n_1-k} &   &  &  \\
          & I_k &   &  \\
          &    & I_k & Y \\
          &     &   & I_{n_2-k} \end{bmatrix}$ for $Y$ and $X$ with coefficients in $K$, belong to $N_kN$. This is immediate multiplying on the left by respectively 
$\begin{bmatrix}
I_{n_1-k} &   & X^\sigma &  \\
          & I_k &   &  \\
          &    & I_k &  \\
          &     &   & I_{n_2-k} \end{bmatrix}$ and $\begin{bmatrix}
I_{n_1-k} &   &  &  \\
          & I_k &   &  Y^\sigma\\
          &    & I_k &  \\
          &     &   & I_{n_2-k} \end{bmatrix}$. 
\end{proof}

Now if the representation $\Delta_1\times\Delta_2$ is distinguished, denoting $\Delta_1\otimes\Delta_2$ by $\Delta$, then at least one of the factors 
$ind_{u^{-1}Pu\cap H}^H((\delta_P^{1/2}\Delta )^u)$ admits on its space a nonzero $H$-invariant linear form. This is equivalent to say that the representation 
$ind_{P\cap uHu^{-1}}^{uHu^{-1}}(\delta_P^{1/2}\Delta )$ admits on its space a nonzero $uHu^{-1}$-invariant linear form. From Frobenius reciprocity law, the space 
$$Hom_{uHu^{-1}}(ind_{P\cap uHu^{-1}} ^{uHu^{-1}}(\delta_P^{1/2}\Delta),1)$$ is isomorphic as a vector space, to 
$$Hom_{P\cap uHu^{-1}}(\delta_P^{1/2}\Delta,\delta_{P\cap uHu^{-1}})= Hom_{P\cap uHu^{-1}}(
 {\delta_P^{1/2}}/{\delta_{P\cap uHu^{-1}}}\Delta, 1).$$\\
 Hence there is on the space $V_{\Delta}$ of $\Delta$ a linear nonzero form $L$, such that for every $p$ in $P\cap uHu^{-1}$, and for every $v$ in $V_{\Delta }$, 
one has $L(\chi(p)\Delta(p)v)=L(v)$, where $\chi(p)=\frac{\delta_P^{1/2}}{\delta_{P\cap uHu^{-1}}}(p)$. As both $\delta_P^{1/2}$ and $\delta_{P\cap uHu^{-1}}$ are trivial 
on $N_k$, so is $\chi$. Now, fixing $k$ such that $u=u_k$, let $n$ belong to $U_k$, from Corollary \ref{util}, we can write $n$ as a product $n_k n_0$, with $n_k$ in $N_k$, 
and $n_0$ in $N$. As $N$ is included in $Ker(\Delta)$, one has $L(\Delta(n)(v))=L(\Delta(n_k n_0)(v))=L(\Delta(n_k )(v))=L(\chi(n_k)\Delta(n_k )v)= L(v)$. Hence $L$ is actually
 a nonzero linear form on the Jacquet module of $V_{\Delta}$ associated with $U_k$. But we also know that $L(\chi(m_k)\Delta(m_k)v)=L(v)$ for $m_k$ in $M_k(F)$, which
 reads according to Lemma \ref{structure}: $L(\delta_{P_k}^{-1/2}(m_k)\Delta(m_k)v)=L(v)$.\\
This says that the linear form $L$ is $M_k(F)$-distinguished on the normalized Jacquet module $r_{ M_k,M}(\Delta)$ (as $M_k$ is also the standard Levi subgroup of $M$).\\
But from Proposition 9.5 of \cite{Z}, there exist quasi-square-integrable representations $\Delta_1'$ of $G_{n_1-k}(K)$, $\Delta_1''$ and $\Delta_2'$ of $G_{k}(K)$, and 
$\Delta_2''$ of $G_{n_2-k}(K)$, such that $\Delta_1=[\Delta_1',\Delta_1'']$ and $\Delta_2=[\Delta_2',\Delta_2'']$, and the normalized Jacquet module $r_{ M_k,M}(\Delta)$ 
is isomorphic to $\Delta_1'\otimes\Delta_1''\otimes\Delta_2'\otimes\Delta_2''$. This latter representation being distinguished by $M_k(F)$, the representations $\Delta_1'$
 and $\Delta_2''$ are distinguished and we have ${\Delta_2'}^\vee={\Delta_1''}^\sigma$. Now we recall from Proposition 12 of \cite{F2}, that we also know that either
 $\Delta_1$ and $\Delta_2$ are Galois auto dual, or we have $\Delta_2^\vee=\Delta_1^\sigma$. In the first case, the representations $\Delta_1$ and $\Delta_2$ are unitary 
because so is their central character, and if nonzero, $\Delta_1'$ and $\Delta_2''$ are also unitary. This implies that either $\Delta_1=\Delta_1'$ and
 $\Delta_2=\Delta_2''$ (i.e. $\Delta_1$ and $\Delta_2$ distinguished), or $\Delta_1=\Delta_1''$ and $\Delta_2=\Delta_2''$ (i.e. $\Delta_1^\sigma=\Delta_2^\vee$).
 This ends the proof of Theorem \ref{cndist}.
\end{proof}

We refer to Section 2 of \cite{M3} for a survey about Asai $L$-functions of generic representations, we will use the same notations here. We recall that if $\pi$ is a generic 
representation of $G_{r}(K)$ for some positive integer $r$, its Asai $L$-function is equal to the product $L_{F,{rad(ex)}}^K(\pi) L_{F,(0)}^K(\pi)$, where 
$L_{F,rad(ex)}^K(\pi)$ is the Euler factor with simple poles, which are the $s_i$'s in $\mathbb{C}/(\frac{2i\pi}{ln(q_F)}\mathbb{Z})$ such that $\pi$ is 
$|\ |_F^{-s_i}$-distinguished, i.e. the exceptional poles of the Asai $L$-function $L_F^K(\pi)$. We denote by $L_{F,ex}^K(\pi)$ the exceptional part of $L_F^K(\pi)$, i.e. the 
Euler factor whose poles are the exceptional poles of $L_F^K(\pi)$, occurring with order equal the order of their occurrence in  $L_F^K(\pi)$. If $\pi'$ is another generic 
representation of $G_{r}(K)$, we denote by $L_{rad(ex)}(\pi\times \pi')$ the Euler product with simple poles, which are the exceptional poles of $L(\pi\times \pi')$ 
(see \cite{CP}, 3.2. Definition). An easy consequence of this definition is the equality $L(\pi\times \pi')=L_{(0)}(\pi\times \pi')L_{rad(ex)}(\pi\times \pi')$. A pole 
$s_0$ of $L(\pi\times \pi')$ is exceptional if and only $\pi'^\vee=|\ |_K^{s_0}\pi$, though only the implication ($s_0$ exceptional $\Rightarrow$ $\pi'^\vee=|\ |_K^{s_0}\pi$)
 is proved in \cite{CP}, the other implication follows from a straightforward adaptation of Theorem 2.2 of \cite{M2}, using Theorem A of \cite{Ber}, instead of using 
Proposition 1.1 (which is actually Ok's theorem) of \cite{M2}.\\
We refer to Definition 3.10 of \cite{M3} for the definition of general position, and recall from Definition-Proposition of \cite{M3}, that if $\Delta_1$ and $\Delta_2$ are two
 square integrable representations of $G_{n_1}(K)$ and $G_{n_2}(K)$, the representation $\Delta_1 |.|_K^{u_1}\times\Delta_2 |.|_K^{u_2}$ is in general position outside a finite
 number of hyperplanes of $(\frac{\mathbb{C}}{2i\pi/Ln(q_F)\mathbb{Z}})^2$ in $(u_1,u_2)$.\\
We refer to Proposition 2.3 of \cite{AKT} and the discussion preceding it for a summary about Bernstein-Zelevinsky derivatives. We use the same notations, except that we use 
the notation $[\rho|\ |_K^{l-1},\dots,\rho]$ where they use the notation $[\rho,\dots,\rho|\ |_K^{l-1}]$.\\
According to Theorem 3.6 of \cite{M3}, we have:

\begin{prop}\label{fact}
Let $m$ be a positive integer, and $\pi$ be a generic representation of $G_m(K)$ such that its derivatives are completely reducible, the Euler factor  $L_{F,(0)}^K(\pi)$ 
(resp. $L_{F}^K(\pi)$) is equal to the l.c.m. $\vee_{k,i}L_{F,ex}^K(\pi_i^{(k)})$ taken over $k$ in $\{1,\dots,n\}$ (resp. in $\{0,\dots,n\}$) and $\pi_i^{(k)}$ in the 
irreducible components of $\pi^{(k)}$.
\end{prop}

An immediate consequence is: 
\begin{cor}\label{factrad}
Let $m$ be a positive integer, and $\pi$ be a generic representation of $G_m(K)$ such that its derivatives are completely reducible, the Euler factor $L_{F,(0)}^K(\pi)$ 
(resp. $L_{F}^K(\pi)$) is equal to the l.c.m. $\vee_{k,i}L_{F,rad(ex)}^K(\pi_i^{(k)})$ taken over $k$ in $\{1,\dots,n\}$ (resp. in $\{0,\dots,n\}$) and $\pi_i^{(k)}$ in the 
irreducible components of $\pi^{(k)}$. 
\end{cor}
\begin{proof}
Let $s$ be a pole of $L_{F,ex}^K(\pi_{i_0}^{(k_0)})$ for $k_0$ in $\{1,\dots,n\}$ and $\pi_{i_0}^{(k)}$ a irreducible component of $\pi^{(k_0)}$. Either $s$ is a pole of
 $L_{F,rad(ex)}^K(\pi_{i_0}^{(k_0)})$, or it is a pole of $L_{F,(0)}^K(\pi_{i_0}^{(k_0)})$, which from Proposition \ref{fact}, implies that it is a pole of some function
 $L_{F,ex}^K((\pi_j^{(k')})$, for $k'>k_0$ and $\pi_j^{(k')}$ a irreducible component of $\pi^{(k')}$. Hence in the factorization 
$L_{F,(0)}^K(\pi)=\vee_{k,i}L_{F,ex}^K(\pi_i^{(k)})$, the factor $L_{F,ex}^K(\pi_{i_0}^{(k_0)})$ can be replaced by $L_{F,rad(ex)}^K(\pi_{i_0}^{(k_0)})$, and the conclusion 
follows from a repetition of this argument. The case of $L_{F}^K(\pi)$ is similar.
\end{proof}

This corollary has a split version:
\begin{prop}\label{factrad2}
Let $m$ be a positive integer, and $\pi$ and $\pi'$ be two generic representations of $G_m(K)$ such that their derivatives are completely reducible, the Euler factor
 $L_{F,(0)}^K(\pi\times\pi')$ (resp. $L_{F}^K(\pi\times\pi')$) are equal to the l.c.m. $\vee_{k,i,j}L_{F,rad(ex)}^K(\pi_i^{(k)}\times{\pi'}_j^{(k)})$ taken over $k$ in 
$\{1,\dots,n\}$ (resp. in $\{0,\dots,n\}$), $\pi_i^{(k)}$ in the irreducible components of $\pi^{(k)}$, and ${\pi'}_j^{(k)}$ in the irreducible components of $\pi^{(k)}$. 
\end{prop}
\begin{proof}
It follows the analysis preceding Proposition 3.3 of \cite{CP}, that one has the equality $L_{(0)}(\pi\times\pi')=\vee_{k,l,i,j}L_{ex}^K(\pi_i^{(k)}\times{\pi'}_j^{(k)})$, 
and the expected statement is a consequence of the argument used in the proof of Corollary \ref{factrad}.

\end{proof}

If $\pi$ is a representation of $G_m(K)$ for some positive integer $m$, admitting a central character, we denote by $R(\Pi)$ the finite subgroup of elements $s$ in 
$\mathbb{C}/(2i\pi/Ln(q_K)\mathbb{Z})$ such that $\pi|\ |_K^s$ is isomorphic to $\pi$.\\
  A consequence of Corollary \ref{factrad} and Theorem \ref{cndist} is the following proposition:

\begin{prop}\label{key}
Let $\Delta$ be a square-integrable representation of $G_n(K)$, and $t$ be a complex number of real part near $n$, 
then the Euler factor $L_F^K(\Pi_t,s)$ equals $L_F^K(\Delta^\sigma,s+2t)L_F^K(\Delta^\vee,s-2t)L(\Delta^\sigma \times\Delta^{\sigma\vee},s)$, and the Euler factor 
$L_{F,(0)}^K(\Pi_t,s)$ equals $\prod_{s_i \in R(\Delta)}( 1-q_K^{s_i-s})L_F^K(\Delta^\sigma,s+2t)L_F^K(\Delta^\vee,s-2t)L(\Delta^\sigma \times\Delta^{\sigma\vee},s)$.
\end{prop}
\begin{proof} 

We first show that for $t$ near $n$, the representation $\Pi_t$ is in general position. According to Definition 4.13 of \cite{M3}, since for such a $t$, $\Pi_t$ is
irreducible, the representation will be in general position if for each $k$ between $1$ and $2n$, the central characters of the irreducible sub quotients of 
$\Pi_t^{(k)}$ have different central characters, and if for each $i$ and $j$ between $1$ and $2n$, the function 
$L((\Delta^\sigma|\ |_K^t)^{(i)}\times (\Delta^{\sigma\vee}|\ |_K^{-t})^{(j)},s)$ has a pole in common, neither with $L_F^K((\Delta^\sigma|\ |_K^t)^{(i)},s)$, nor
 with $L_F^K((\Delta^\sigma|\ |_K^{-t})^{(j)},s)$. According to Corollary 4.4 and Remark 4.5 of \cite{M3}, the latter condition is equivalent to the fact that the function 
$L( \Delta^\sigma|\ |_K^t \times  \Delta^{\sigma\vee}|\ |_K^{-t} ,s)$ has a pole in common, neither with $L_F^K(\Delta^\sigma|\ |_K^t ,s)$, nor
 with $L_F^K(\Delta^{\sigma\vee}|\ |_K^{-t},s)$, and by invariance of the L-functions under $\sigma$, we can remove it in the preceding expressions.\\
We start by proving the assumption on the central characters. Writing the discrete series representation 
$\Delta$ under the form $St_{l}(\rho)=[\rho|\ |_K^{(l-1)/2},\dots,\rho|\ |_K^{(1-l)/2}]$ for a positive integer $l$ and a unitary supercuspidal representation 
$\rho$ of $G_m(K)$, with $lm=n$, from Proposition 9.6 of \cite{Z}, the derivative $\Pi_t^{(k)}$ is zero unless $k$ is of the form $mk'$ for $k'$ between $1$ and $l$,
 in which case its irreducible components are the 
$$|\ |_K^t [\rho^\sigma|\ |_K^{(l-1)/2},\dots,\rho^\sigma|\ |_K^{(1-l)/2 +i}]\times |\ |_K^{-t} [\rho^\vee|\ |_K^{(l-1)/2},\dots,\rho^\vee|\ |_K^{(1-l)/2 +k'-i}]$$
for $i$ between $0$ and $k'$. The exponent of the central character of this representation is equal to $$m[Re(t)(k'-2i)+(l-i)i/2 + (l+i-k')(k'-i)/2].$$
One the checks that for $i'\neq i$, the two exponents are different for $Re(t)$ near $n$.\\
Concerning the condition on the L functions, it follows from the proof of Proposition 4.16 in \cite{M3}, that if
 $L( \Delta |\ |_K^t \times  \Delta^{ \vee}|\ |_K^{-t} ,s)$ has a pole in common with $L_F^K(\Delta |\ |_K^t ,s)$, then one would have 
$\rho^{\sigma \vee}=\rho|\ |_K^{a+2t}$ for $a$ an integer between $-n$ and $n$, which is impossible because $\rho$ has a unitary central character and $t$ is near $n$.
 We obtain a similar contradiction if we assume that $L( \Delta |\ |_K^t \times  \Delta^{ \vee}|\ |_K^{-t} ,s)$ has a pole in common with 
$L_F^K(\Delta^{\sigma\vee}|\ |_K^t,s)$.\\
Hence for $t$ near $n$, the representation $\Pi_t$ is in general position.\\

Now from Corollary \ref{factrad}, we know that given the hypothesis of the proposition, the function $L_F^K(\Pi_t,s)$ is equal to the l.c.m. 
$\vee_{k_1, k_2/k_1+k_2\geq1}L_{F,rad(ex)}^K((\Delta^\sigma|\ |_K^t) ^{(k_1)}\times (\Delta^\vee|\ |_K^{-t})^{(k_2)})$. Writing the discrete series representation 
$\Delta$ under the form $St_{l}(\rho)=[\rho|\ |_K^{(l-1)/2},\dots,\rho|\ |_K^{(1-l)/2}]$ for a positive integer $l$ and a unitary supercuspidal representation 
$\rho$ of $G_m(K)$, with $lm=n$, the representation $(\Delta^\sigma|\ |_K^t) ^{(k_1)}$ (resp. $(\Delta^\vee|\ |_K^{-t})^{(k_2)})$) is equal to zero unless there exists
 an integer $k'_1$ with $k_1=mk'_1$ (resp. $k'_2$ with $k_2=mk'_2$), in which case it is equal to $St_{l-k'_1}(\rho^\sigma)|\ |_K^{k'_1/2+t}$ 
(resp. $St_{l-k'_2}(\rho^\vee)|\ |_K^{k'_2/2-t}$).\\
Suppose that the representations $(\Delta^\sigma|\ |_K^t) ^{(k_1)}$ and $(\Delta^\vee|\ |_K^{-t})^{(k_2)})$ are not zero (hence $k_i=mk'_i$ for a integer $k'_i$), 
a complex number $s_0$ is a pole of $L_{F,rad(ex)}^K((\Delta^\sigma|\ |_K^t) ^{(k_1)}\times (\Delta^\vee|\ |_K^{-t})^{(k_2)})$ if and only if the representation 
$St_{l-k'_1}(\rho^\sigma)|\ |_K^{k'_1/2+t}\times St_{l-k'_2}(\rho^\vee)|\ |_K^{k'_2/2-t}$ is $|\ |_K^{-s_0}$-distinguished, i.e. 
$St_{l-k'_1}(\rho^\sigma)|\ |_K^{(k'_1+s_0)/2+t}\times St_{l-k'_2}(\rho^\vee)|\ |_K^{(k'_2+s_0)/2-t}$ is distinguished. But from Theorem \ref{cndist}, this implies that 
$k'_1$ and $k'_2$ are equal to an integer $k'$, (i.e. $k_1=k_2=k$), and that the image of $s_0+k'$ in $\mathbb{C}/(2i\pi/Ln(q_K)\mathbb{Z})$ belongs to the group 
$R(St_{l-k'}(\rho))=R(\rho)$ (in particular, we have $Re(s_0+k')=0$).\\
 Conversely if this is the case, the representation 
$St_{l-k'}(\rho^\sigma)|\ |_K^{(k'+s_0)/2+t}\times St_{l-k'}(\rho^\vee)|\ |_K^{(k'+s_0)/2-t}$ which is equal to 
$St_{l-k'}(\rho^\sigma)|\ |_K^{(k'+s_0)/2+t}\times St_{l-k'}(\rho^\vee)|\ |_K^{(-k'-s_0)/2-t}$, is distinguished from Theorem \ref{dist1}, as 
$Re((k'+s_0)/2+t)=Re(t)>n/2\geq (n-k)/2$ for $Re(t)$ near $n$.\\

 Hence nontrivial Euler factors $L_{F,rad(ex)}^K((\Delta^\sigma|\ |_K^t) ^{(k_1)}\times (\Delta^\vee|\ |_K^{-t})^{(k_2)})$ belong to one of the three following classes:

\begin{enumerate}
 \item $L_{F,rad(ex)}^K((\Delta^\sigma|\ |_K^t) ^{(k_1)})$ for $k_2=n$ and $k_1\geq 0$. In this case, if $L_{F,rad(ex)}^K((\Delta^\sigma|\ |_K^t) ^{(k_1)})$ is not $1$, 
it is equal to $L_{F,rad(ex)}^K(St_{l-k'_1}(\rho^\sigma)|\ |_K^{k'_1/2+t})$ for $k_1=mk'_1$, and a pole $s_0$ of this function is such that 
$St_{l-k'_1}(\rho^\sigma)|\ |_K^{(s_0+k'_1)/2+t}$ is distinguished, hence considering central characters, we have $Re_(s_0)=-k'_1-2Re(t)<-n$.

\item $L_{F,rad(ex)}^K((\Delta^\vee|\ |_K^{-t}) ^{(k_2)})$ for $k_1=n$ and $k_2\geq 0$. In this case, if $L_{F,rad(ex)}^K((\Delta^\sigma|\ |_K^t) ^{(k_2)})$ is not $1$, 
it is equal to $L_{F,rad(ex)}^K(St_{l-k'_2}(\rho^\vee)|\ |_K^{k'_2/2-t})$ for $k_2=mk'_2$, and a pole $s_0$ of this function is such that 
$St_{l-k'_2}(\rho^\sigma)|\ |_K^{(s_0+k'_2)/2-t}$ is distinguished, hence considering central characters, we have $Re_(s_0)=-k'_2+2Re(t)>0$.

\item $L_{F,rad(ex)}^K((\Delta^\sigma|\ |_K^t) ^{(k_3)}\times (\Delta^\sigma|\ |_K^{-t}) ^{(k_3)})$ for $k_1=k_2=k_3\geq 1$. In this case, if the Euler factor is not $1$, 
we know that we have $Re(s_0)=-k'_3$ for $k'_3$ in $\{0,\dots,n/m\}$ verifying $k_3=mk'_3$, or more precisely that the image of $s_0+k'_3$ in 
$\mathbb{C}/(2i\pi/Ln(q_K)\mathbb{Z})$ belongs to the group $R(St_{l-k'_3}(\rho))=R(St_{l-k'_3}(\rho^\sigma))$. This is equivalent to the relation 
$[{\Delta^{\sigma\vee}} ^{(k_3)}]^\vee= |\ |_K^{s_0}(\Delta^\sigma)^{(k_3)}$, which is itself equivalent to the fact that $s_0$ is a pole of 
$L_{rad(ex)}({\Delta^{\sigma}}^{(k_3)}\times {\Delta^{\sigma\vee}}^{(k_3)})$ (see Th. 1.14 of \cite{M3}), hence we have 
$L_{F,rad(ex)}^K((\Delta^\sigma|\ |_K^t) ^{(k_3)}\times (\Delta^\sigma|\ |_K^{-t}) ^{(k_3)})=L_{rad(ex)}({\Delta^{\sigma}}^{(k_3)}\times {\Delta^{\sigma\vee}}^{(k_3)})$. 

\end{enumerate}

In particular, two non trivial factors that don't belong to the same class have no pole in common. We deduce that the Euler factor $L_{F,(0)}^K(\Pi_t,s)$ 
is equal to 
$$[\vee_{k_1}L_{F,rad(ex)}^K((\Delta^\sigma|\ |_K^t) ^{(k_1)}][\vee_{k_2}L_{F,rad(ex)}^K((\Delta^\vee|\ |_K^{-t}) ^{(k_2)}]
[\vee_{k_3} L_{rad(ex)}({\Delta^{\sigma}}^{(k_3)}\times {\Delta^{\sigma\vee}}^{(k_3)})]$$
 for $k_1\geq 0$, $k_2\geq 0$ and $k_3\geq 1$. The two first factors are respectively equal to $L_F^K(\Delta^\sigma|\ |_K^t)$ and 
$ L_F^K(\Delta^\vee|\ |_K^{-t})$ according to Corollary \ref{factrad}, and the third factor is equal from Proposition \ref{factrad2} to 
$L_{(0)}(\Delta^{\sigma}\times \Delta^{\sigma\vee})$, which is itself equal to 
$L(\Delta^{\sigma}\times \Delta^{\sigma\vee})/L_{rad(ex)}(\Delta^{\sigma}\times \Delta^{\sigma\vee})$. We then notice that $s_0$ is an
 exceptional pole of $L(\Delta^{\sigma}\times \Delta^{\sigma\vee})$ if and only if its image in $\mathbb{C}/(2i\pi/Ln(q_K)\mathbb{Z})$ belongs to 
$R(\Delta)$, which implies the equality $L_{rad(ex)}(\Delta^{\sigma}\times \Delta^{\sigma\vee})=1/\prod_{s_i \in R(\Delta)}( 1-q^{s_i-s})$. Hence we 
deduce the equalities 
$$\begin{array}{ll} L_{F,(0)}^K(\Pi_t,s) &=L_F^K(\Delta^\sigma|\ |_K^t,s)L_F^K(\Delta^\vee|\ |_K^{-t},s)
[L(\Delta^{\sigma}\times \Delta^{\sigma\vee},s)/L_{rad(ex)}(\Delta^{\sigma}\times \Delta^{\sigma\vee},s)]\\   
 &=\prod_{s_i \in R(\Delta)} (1-q_K^{s_i-s})L_F^K(\Delta^\sigma|\ |_K^t,s)L_F^K(\Delta^\vee|\ |_K^{-t},s)L(\Delta^{\sigma}\times \Delta^{\sigma\vee},s)\end{array}$$
The second statement of the proposition follows, as tensoring by $|\ |^u$ the representation, is equivalent to make a translation by $2u$ of the 
Asai $L$ function.\\
As the function $L_F^K(\Pi_t,s)$ is equal to the product $L_{F,rad(ex)}^K(\Pi_t,s)L_{F,(0)}^K(\Pi_t,s)$. It remains to show that the function 
$L_{F,rad(ex)}^K(\Pi_t,s)$ is equal to the factor $\prod_{s_i \in R(\Delta)}1/( 1-q_K^{s_i-s})$. But we already know that it is equal to the product of 
the $1/(1-q^{s_i-s})$'s, for $s_i$'s such that $\Pi_t$ is $|\ |_F^{-s_i}$-distinguished. As $\Pi_t$ is $|\ |_F^{-s_i}$-distinguished if and only if 
$\Pi_t |\ |_K^{s_i/2}=\Delta^{\sigma}|\ |_K^{t+s_i/2}\times \Delta^{\vee}|\ |_K^{-t+s_i/2}$ is distinguished, Theorem \ref{cndist} implies that 
if $\Pi_t$ is $|\ |_F^{-s_i}$-distinguished, either we 
have $\Delta^{\sigma}|\ |_K^{t+s_i/2}$ and $\Delta^{\vee}|\ |_K^{-t+s_i/2}$ distinguished (hence Galois-auto-dual), or we have 
$(\Delta^{\sigma}|\ |_K^{t+s_i/2})^\sigma= (\Delta^{\vee}|\ |_K^{-t+s_i/2})^\vee$. The first case cannot occur because quasi-square-integrable 
distinguished representations must be unitary (because distinguished representations have unitary central character), and this would imply 
$Re(t+s_i/2)=Re(t-s_i/2)=0$, which would in turn imply $Re(t)=0$. The second case clearly implies that $s_i$ belongs to $R(\Delta)$. Conversely, if 
$s_i$ belongs to $R(\Delta)$, its real part is zero, and it is immediate that the representation $\Pi_t|\ |_K^{s_i/2}$ verifies the hypothesis of Theorem 
\ref{dist1}. This concludes the proof of the first statement.

\end{proof}

\begin{df-prop}\label{simple}
We denote by $P_{(0)}(\Pi,t,s)$ the element of $\mathbb{C}[q_F^{\pm t},q_F^{\pm s}]$ defined by the expression 
$\frac{\prod_{s_i \in R(\Delta)}(1-q^{s_i-s})}{L_F^K(\Delta^\sigma,s+2t)L_F^K(\Delta^\vee,s-2t)L(\Delta^\sigma \times\Delta^\vee,s)}$. 
The expression $P_{(0)}(\Pi,t,1)$ defines a nonzero element of $\mathbb{C}[q_F^{\pm t}]$, having simple roots. 
For any complex number $t_0$, the expression $P_{(0)}(\Pi,t_0,s)$ defines a nonzero element of $\mathbb{C}[q_F^{\pm t}]$, having an 
at most simple root at $s=1$.
\end{df-prop}
\begin{proof} 
As the $s_i$'s have real part equal to zero, and as the function $L(\Delta^\sigma \times\Delta^{\sigma\vee},s)$ admits no pole for $Re(s)>0$ 
(see \cite{JPS1}, 8.2 (6)), the constant $c=\frac{\prod_{s_i \in R(\Delta)}(1-q^{s_i-1})}{L(\Delta^\sigma \times\Delta^{\sigma\vee},1)}$ is nonzero.
 Hence the zeros of $P_{(0)}(\Pi,t,1)$ are the poles of $L_F^K(\Delta^\sigma,1+2t)L_F^K(\Delta^\vee,1-2t)$. From Proposition 3.1 of \cite{M3}, the 
function $L_F^K(\Delta^\sigma,1+2t)$ has simple poles which occur in the domain $Re(1+2t)<0$ whereas the function $L_F^K(\Delta^\vee,1-2t)$ has simple 
poles which occur in the domain $Re(1-2t)<0$, hence those two functions have no common pole, and there product have simple poles. The second part is a
 consequence of the fact that the function $L_F^K(\Delta^\sigma,s+2t_0)$ has simple poles, and if it has a pole at $1$, then $Re(1+2t_0)<0$, whereas 
$L_F^K(\Delta^\vee,s-2t_0)$ also has simple poles, and if it has a pole at $1$, then $Re(1-2t_0)<0$, so that both cannot have a pole at $1$ at the same 
time.
\end{proof}

\begin{LM}\label{poles}
For every $f$ in $\mathcal{F}_{\Pi}$, the expression $P_{(0)}(\Pi,t,s) I_{(0)}(W_{f_t},s)$ defines an element of $\mathbb{C}[q_F^{\pm t},q_F^{\pm s}]$. 
This implies that for fixed $f$ in $\mathcal{F}_{\Pi}$, the function  $I_{(0)}(W_{f_t},1)$ is well defined and belongs to $\mathbb{C}(q_F^{t})$, and for $t_0$ in 
$\mathbb{C}$, the function $I_{(0)}(W_{f_{t_0}},s)$ is well defined and belongs to $\mathbb{C}(q_F^{-s})$. Moreover the function $I_{(0)}(W_{f_t},1)$ has a pole 
at $t_0$ in $\mathbb{C}$, if and only if the function $I_{(0)}(W_{f_{t_0}},s)$ in $\mathbb{C}(q_F^{-s})$ has a pole at $1$, in which case the couple $(t_0,1)$ 
lies in a polar locus of the function $P_{(0)}(\Pi,t,s)$. In this case the functions $P_{(0)}(\Pi,t,1) I_{(0)}(W_{f_{t}},1)$ and 
$P_{(0)}(\Pi,t_0,s) I_{(0)}(W_{f_{t_0}},s)$ have the same limit when $t$ tends to $t_0$ and $s$ tends to $1$, which is nonzero. 
\end{LM}

\begin{proof}
Let $f$ be in $\mathcal{F}_{\Pi}$, the function $P_{(0)}(\Pi,t,s)I_{(0)}(W_{f_t},s)$ belongs to $\mathbb{C}(q_F^{-t}, q_F^{-s})$, hence it is the quotient
 of two polynomials $P(q_F^{-t}, q_F^{-s})/Q(q_F^{-t}, q_F^{-s})$. If $Q$ is not constant, we write $Q(q_F^{-t}, q_F^{-s})$ under the form 
$\sum_{i\in I} a_i(q_F^{-t})q_F^{-is}$, with $I$ a finite subset of $\mathbb{Z}$, and the $a_i$'s in $\mathbb{C}[X]-{\{0\}}$. There are two real numbers $\alpha<\alpha'$  
such that $[\alpha,\alpha']$ is a subset of a neighbourhood of $n$ containing no real part of a zero of the function $t\mapsto a_{i_0}(q_F^{-t})$, for $i_0$ the minimum of $i$.
As the functions $a_{i}(q_F^{-t})$ are bounded for $Re(t)\in [\alpha,\alpha']$, there is a real number $r$, such that for $Re(t)\in [\alpha,\alpha']$, and $Re(s)\geq r$, 
the function 
$P_{(0)}(\Pi,t,s)I_{(0)}(W_{f_t},s)$ is given by an absolutely convergent Laurent development $\sum_{k\geq n_{0}} c_k(t) q_F^{-ks}$ with $c_k$ in 
$\mathbb{C}[q_F^{\pm t}]$. Moreover if we choose $[\alpha,\alpha']$ so that $\Pi_t$ satisfies the hypothesis of Proposition \ref{key} for $Re(t)$ in $[\alpha,\alpha']$, 
then for fixed $t$ with $Re(t)$ in $[\alpha,\alpha']$, the function 
$P_{(0)}(\Pi,t,s) I_{(0)}(W_{f_t},s)=I_{(0)}(W_{f_t},s)/L_{(0)}(\Pi_t,s)$ actually belongs to $\mathbb{C}[q_F^{\pm s}]$. Suppose there were an infinite 
number of nonzero $c_k$'s, then for $t$ of real part in $[\alpha,\alpha']$, and outside the countable number of zeros of the $c_k$'s, and $Re(s)$ large,
 the Laurent development 
$\sum_{k\leq n_{0}} c_k(t) q_F^{-ks}$ would not be finite, a contradiction. Hence for $f$ in $\mathcal{F}_{\Pi}$, the function 
$P_{(0)}(\Pi,t,s)I_{(0)}(W_{f_t},s)$ defines an element of $\mathbb{C}[q_F^{\pm t},q_F^{\pm s}]$. \\
Now the function $I_{(0)}(W_{f_t},1)$ defines an element of $\mathbb{C}(q_F^{-t})$ whose poles form a subset of the poles of $1/P_{(0)}(\Pi,t,1)$, and for 
$t_0$ in $\mathbb{C}$, the function $I_{(0)}(W_{f_{t_0}},s)$ defines an element of $\mathbb{C}(q_F^{-s})$ whose poles form a subset of the poles of 
$1/P_{(0)}(\Pi,t_0,s)$.\\
For the final statement, if $t_0$ is a pole of $I_{(0)}(W_{f_{t}},1)$, then it must be a zero of the function $P_{(0)}(\Pi,t,1)$, which is simple 
according to Definition-Proposition \ref{simple}, as $P_{(0)}(\Pi,t,1)I_{(0)}(W_{f_{t}},1)$ is polynomial, the pole $t=t_0$ is also simple. Hence the function 
$P_{(0)}(\Pi,t,1)I_{(0)}(W_{f_{t}},1)$ has nonzero limit when $t$ tends to $t_0$.
 As the function $P_{(0)}(\Pi,t,s) I_{(0)}(W_{f_t},s)$ belongs to $\mathbb{C}[q_F^{\pm t},q_F^{\pm s}]$, the function $P_{(0)}(\Pi,t_0,s) I_{(0)}(W_{f_{t_0}},s)$ 
tends to the same limit when $s$ tends to $1$. Conversely if $1$ is a pole of  $I_{(0)}(W_{f_{t_0}},s)$, then it must be a zero of the function 
$P_{(0)}(\Pi,t_0,s)$, which is simple according to Definition-Proposition \ref{simple}, as 
$P_{(0)}(\Pi,t_0,s)I_{(0)}(W_{f_{t_0}},s)$ is polynomial, the pole $s=1$ is also simple. Hence the function $P_{(0)}(\Pi,t_0,s)I_{(0)}(W_{f_{t_0}},s)$ has 
nonzero limit when $s$ tends to $1$. As the function $P_{(0)}(\Pi,t,s) I_{(0)}(W_{f_t},s)$ belongs to $\mathbb{C}[q_F^{\pm t},q_F^{\pm s}]$, the 
function $P_{(0)}(\Pi,t,1) I_{(0)}(W_{f_{t}},1)$ tends to the same limit when $t$ tends to $t_0$.
\end{proof}

Finally we can prove the main result.

\begin{thm}\label{dist2}
Let $\Delta'$ be a quasi-square-integrable representation of $G_n(K)$, then the representation $\Delta'^\sigma\times\Delta'^\vee$ of $G_{2n}(K)$ is 
distinguished. 
\end{thm}

\begin{proof}
Write $\Delta'=\Delta |.|_K^u$, for $\Delta$ a square-integrable representation, and $u$ a complex number.
Denoting by $\Pi_t$ the representation $\Delta^\sigma|.|_K^t\times\Delta^\vee|.|_K^{-t}$, we know from Proposition \ref{dist1} that $\Pi_t$ is 
distinguished for $Re(t)$ near $n$. Hence for $Re(t)$ near $n$, we know from Proposition \ref{invform}, that the linear form 
$W_{f_t} \mapsto \underset{ s \rightarrow 1}{lim}  I_{(0)}(W_{f_{t}},s)/L_{(0)}(\Pi_t,s)$ is nonzero and $G_{2n}(F)$-invariant.\\
Suppose that $t$ is in a neighbourhood of $n$ such that $\Pi_t$ is in general position (the see proof of Proposition \ref{key}), then the function 
$1/L_{(0)}(\Pi_t,s)$ is equal to $P_{(0)}(\Pi,t,s)$. But the function $P_{(0)}(\Pi,t,1)$, which is a nonzero polynomial in $q_F^{-t}$, has no zeros for $Re(t)$ 
in some open subset in a neighbourhood of $n$. From this we deduce that for $Re(t)$ in this open subset, according to Lemma \ref{poles}, 
the functions $s\mapsto I_{(0)}(W_{f_{t}},s)$ and 
$t'\mapsto  I_{(0)}(W_{f_{t'}},1)$ have respectively no pole at $s=1$ and $t'=t$, and we have 
$\underset{ s \rightarrow 1}{lim}  I_{(0)}(W_{f_{t}},s)= \underset{ t' \rightarrow t}{lim}  I_{(0)}(W_{f_{t'}},1)$. Hence for $Re(t)$ in this open subset, 
if $h$ belongs to $G_{2n}(F)$ the two functions $I_{(0)}(W_{f_{t}},1)$ and 
$I_{(0)}(\rho_t(h)W_{f_{t}},1)$ coincide, but as they are rational functions in $q_F^{-t}$, they are equal. Hence for $f$ in the space of $\Pi_0$, and 
$h$ in $G_{2n}(F)$, the functions $I_{(0)}(W_{f_{t}},1)$ and $I_{(0)}(\rho_t(h)W_{f_{t}},1)$ are equal.\\
 Suppose that for every $f$ in the space of $\Pi_0$, the function $I_{(0)}(\rho_t(h)W_{f_{t}},1)$ has no pole at $t=u$, then according to Proposition 
\ref{poles}, for every $f$ in the space of $\Pi_0$, the function $I_{(0)}(\rho_u(h)W_{f_{u}},s)$ has no pole at $s=1$, and if $h$ is in $G_{2n}(F)$, one 
has $\underset{ s \rightarrow 1}{lim}I_{(0)}( \rho_u(h) W_{f_{u}},s)= \underset{ t \rightarrow u}{lim}  I_{(0)}(\rho_t(h)W_{f_{t}},1)=
 \underset{ t \rightarrow u}{lim} I_{(0)}(W_{f_{t}},1)= \underset{ s \rightarrow 1}{lim} I_{(0)}(W_{f_{u}},s)$. 
Hence we have a $G_{2n}(F)$-invariant linear form $f_u\mapsto \underset{ s \rightarrow 1}{lim} I_{(0)}(W_{f_{u}},s)$ on the space of $\Pi_u$. Moreover, as 
$W_{f_{u}}$ describes the space $W(\pi_u,\psi)$ when $f_u$ describes the space of $\Pi_u$, and as the restrictions to $P_n(K)$ of functions of 
$W(\pi_u,\psi)$ form a vector space with subspace $C_c^\infty(N_n(K)\backslash P_n(K), \psi)$, if we choose $W_{f_{u}}$ with restriction to $P_n(K)$ 
positive and in $C_c^\infty(N_n(K)\backslash P_n(K), \psi)$, then we have $I_{(0)}(W_{f_{u}},1)=\int_{N_n(F)\backslash P_n(F)} W_{f_{u}}(p) dp>0$, and 
the $G_{2n}(F)$-invariant linear form defined above is nonzero, hence $\Pi_u=\Delta'^\sigma\times\Delta'^\vee$ is distinguished.\\
 Now if for some $f$ in in the space of $\Pi_0$, the function $I_{(0)}(\rho_t(h)W_{f_{u}},s)$ has a pole at $s=1$, it is a consequence of Lemma 
\ref{poles} that we have $\underset{ s \rightarrow 1}{lim} P_{(0)}(\Pi,u,s)I_{(0)}(W_{f_{u}},s)$ is nonzero, and from the same Lemma, we know that for 
every $f$ in in the space of $\Pi_0$, and $h$ in $G_{2n}(F)$, we have 
$$\begin{aligned} \underset{ s \rightarrow 1}{lim}P_{(0)}(\Pi,u,s)I_{(0)}( \rho_u(h) W_{f_{u}},s)= 
\underset{ t \rightarrow u}{lim} P_{(0)}(\Pi,t,1) I_{(0)}(\rho_t(h)W_{f_{t}},1)\\ 
= \underset{ t \rightarrow u}{lim}  P_{(0)}(\Pi,t,1) I_{(0)}(W_{f_{t}},1)= \underset{ s \rightarrow 1}{lim} P_{(0)}(\Pi,u,s)I_{(0)}(W_{f_{u}},s).\end {aligned}$$  
Hence in this case too, the representation $\Pi_u=\Delta'^\sigma\times\Delta'^\vee$ is distinguished.\end{proof}

\section{Appendix}

In this appendix, we prove a result that seems to be well-known (e.g. it is used in the proof of the lemma of Section 4. in \cite{F2}). 
However we couldn't find a proof in literature. To do this we will give in Proposition \ref{zerobehaviour} and Theorem \ref{restorus} a very precise description of 
the restriction of Whittaker functions on $G_{n-1}(K)$
 to the standard maximal torus of $G_{n-1}(K)$, refining Proposition 2.2 of \cite{JPS1} and Proposition 2.6 of \cite{CP}. As for the proposition, these facts seem 
to be well-known, but there seems to be a lack of references.\\

\begin{prop}\label{l2}
 Let $\pi$ be a smooth submodule of finite length of $C^\infty_c(N_n(K)\backslash P_n(K),\psi)$ such that all of the central 
exponents of its 
shifted derivatives (see \cite{Ber}, 7.2.) are positive, then for any $W$ in $\pi$ the following integral converges:  
$$\int_{N_n(K)\backslash P_n(k)} |W(p)|^2 dp .$$
\end{prop}

Before proving this, we give a precise description of the behaviour near zero of the functions $t\mapsto W \begin{pmatrix}
t I_j & \\  & I_{n-j}                                                                                                                                           
  \end{pmatrix}$, for $t\in K^*$, $j$ between $1$ and $n$, and $W$ in $\pi$. This generalizes Proposition 2.6 of \cite{CP}, where the case of $W$ with 
 $g\mapsto |g|_K^{(n-j-1)/2}W \begin{pmatrix}
g & \\  & I_{n-j}                                                                                                                                           
  \end{pmatrix}$ projecting into an irreducible submodule of $\pi^{(n-j)}$ (or equivalently where the case of completely reducible derivatives) is treated.\\

Let $\Phi^{-}$ be the functor which to a $P_n(K)$-smooth module $(\pi,V)$, associates the $P_{n-1}(K)$-smooth module $$\frac{V}{V(U_n(K),\psi)},$$ where $U_n(K)$ is the 
unipotent radical of $P_{n-1,1} (K)$, and $V(U_n(K),\psi)$ is the vector subspace spanned by the vectors $\pi(u)v-\psi(u)v$ for $v$ in $V$ and $u$ in $U_n(K)$. The action 
of $P_{n-1}(K)$ on $\Phi^{-}(V)$ is its natural action twisted by the character $|\ |_K^{-1/2}$. For $V \subset C^\infty_c(N_n(K)\backslash P_n(K),\psi)$, one sees that
 $V(U_n(K),\psi)$ is the kernel of the restriction map from $C^\infty_c(N_n(K)\backslash P_n(K),\psi)$ to $C^\infty_c(N_{n-1}(K)\backslash P_{n-1}(K),\psi)$.\\
 Applying this repeatedly, one sees (see \cite{CP} Proposition 2.2), one sees that a model of $\pi_{(n-j-1)}=(\Phi^{-})^{n-j-1}(\pi)$, for $\pi$ a 
submodule of $C^\infty_c(N_n(K)\backslash P_n(K),\psi)$, 
is the space of functions of the form $$g\mapsto |g|_K^{-(n-j-1)/2}W \begin{pmatrix}
g & \\  & I_{n-j}                                                                                                                                           
  \end{pmatrix}$$ with $g$ in $G_{j}(K)$, and $W$ in $\pi$, with  $P_{j+1}(K)$ acting by right translation (the twist by $|g|_K^{-(n-j-1)/2}$ appears because of the presence
 of the twist by $|det(\ )|_K^{-1/2}$ in the definition of $\Phi^{-}$, and $g$ is in $G_{j}(K)$ because the quotient space $N_{j+1}(K)\backslash P_{j+1}(K)$ identifies with 
$N_{j}(K)\backslash G_{j}(K)$).\\
Now we introduce the functor $\Psi^{-}$ which to a $P_n(K)$-smooth module $(\pi,V)$, associates the $G_{n-1}(K)$-smooth module $$\frac{V}{V(U_n(K),1)},$$ 
where $V(U_n(K),1)$ is the vector subspace spanned by the vectors $\pi(u)v- v$ for $v$ in $V$ and $u$ in $U_n(K)$.\\
Then one shows (\cite{CP} Proposition 2.3) that when $V$ is the space of $\pi_{(n-j-1)}$ with $\pi$ a submodule of $C^\infty_c(N_n(K)\backslash P_n(K),\psi)$, 
then $V(U_n(K),1)$ is the subspace of 
functions $$g\mapsto |g|_K^{-(n-j-1)/2}W \begin{pmatrix}
g & \\  & I_{n-j}                                                                                                                                           
  \end{pmatrix}$$ which vanish when the last row of $g$ is in a nieghbourhood of zero (depending on $W$).\\
By definition, the $G_j(K)$-module $\pi^{(n-j)}=\Psi^{-}(\pi_{(n-j-1)})=\Psi^{-}(\Phi^{-})^{n-j-1}(\pi) $ is the $(n-j)$-th derivative of $\pi$,
 and the shifted derivative $\pi^{[n-j]}$ is equal to $| \ |_K^{1/2} \pi^{(n-j)}$. 
It is known (\cite{BZ} section 3), that the functor $\Phi^{-}$ and $\Psi^{-}$ are exact and take representations of finite length to representations of finite length.\\

Hence when $\pi$ satisfies the conditions of proposition \ref{l2}, the $G_j(K)$-module $\pi^{(n-j)}$ has finite length. We are going to analyze the 
smooth representation of the center $K^*$ of $G_j(K)$ on the space $E$ of $\tau=\pi^{(n-j)}$. As $E$ is a $G_j(K)$-module $E$ of finite length, it has a filtration 
${0}=E_0 \subset E_1 \subset \dots \subset E_{r-1} \subset E_{r}=E$ with $V_i=E_i/E_{i-1}$ irreducible $G_j(K)$-modules, on which $K^*$ act by the central character $c_i$ 
of the representation of $G_j(K)$. We first show that we can reduce to finite dimension.

\begin{LM}\label{fin-dim}
 Any vector $v$ of $E$ lies in a finite dimensional $K^*$-submodule.
\end{LM}
\begin{proof}
One proves this by induction on the smallest $i$ such that $E_i$ contains $v$. If this $i$ is $1$, the group $K^*$ only multiplies $v$ by a scalar, and we are done.\\
Suppose that the result is known for $E_i$, we take $v$ in $E_{i+1}$ but not in $E_i$, then for every $t$ in $K^*$, the vector $\tau(t)v-c_i(t)v$ belongs to $E_i$. 
By smoothness, the set $\{ \tau(u)v \ | t\in U_K \}$ is actually equal to $\{ \tau(u)v \ | t\in P \}$ for $P$ a finite set of $U_K$. The vector space generated by this set
is stabilized by $U_K$, and has a finite basis $v_1,\dots,v_m$. Now the vector $\tau(\pi_K)v_l-c_i(\pi_K)v_l$ belongs to $E_i$, hence to a finite dimensional 
$K^*$-submodule $V_l$ of $E_i$.
Finally the finite dimensional $Vect(v_1,\dots,v_m)+V_1+\dots+V_m$ is stable under $ \pi_K $ and $ U_K $, hence $K^*$, and contains $v$.
\end{proof}

Each vector of $E$ will thus belong to a subspace satisfying the statement of the following:

\begin{prop}\label{center-repr}
 If $E'$ is a non zero finite dimensional $K^*$-submodule of $E$, then $E'$ has a basis $B$ in which the action of $K^*$ is given by a matrix block diagonal matrix 
$Mat_B (\tau(t))$ with each block of the form:
$$ \begin{pmatrix}
c (t) & c (t)P_{1,2} (v_K(t))     & c (t)P_{1,3} (v_K(t))     & \dots                      & c (t)P_{1,q} (v_K(t))   \\
                    & c (t) & c (t)P_{2,3} (v_K(t))     &  \dots                     & c (t)P_{2,q} (v_K(t))   \\
                    &                    & \ddots              &                            & \vdots            \\  
                    &                    &                     & c (t)      &  c (t)P_{q-1,q} (v_K(t)) \\ 
                    &                    &                     &                            & c(t)                                                                                                     
  \end{pmatrix},$$
for $c$ one of the $c_i$'s, $q$ a positive integer depending on the block, and the $P_{i,j}$'s being polynomials with no constant term of degree at most $j-i$.
\end{prop}

\begin{proof}

First we decompose $E'$ as a direct sum under the action of the compact abelian group $U_K$. Because $E'$ has a filtration by the spaces $E' \cap E_i$,
and that $K^*$ acts on each sub factor as one of the $c_i$'s, the group $U_K$ acts on each weight space as the restriction of one of the $c_i$'s. Now each weight space is
stable under $K^*$ by commutativity, and so we can restrict ourselves to the case where $E'$ is a weight-space of $U_K$.\\
Again $E'$ has a filtration, such that $K^*$ acts on each sub factor as one of the $c_i$'s (with all these characters having the same restriction to $U_K$),
 let's say $c_{i_1}, \dots, c_{i_k}$, in particular, we deduce that the endomorphism $\tau(\pi_K)$ has a triangular matrix in a basis adapted to this filtration, 
with eigenvalues $c_{i_1}(\pi_K), \dots, c_{i_k}(\pi_K)$. As $\tau(\pi_K)$ is trigonalisable, the space $E'$ is the direct sum its characteristic subspaces, and
again these characteristic subspaces are stable under $K^*$.\\
 So finally one can assume that $E'$ is a characteristic subspace for some eigenvalue $c(\pi)$ of $\tau(\pi_K)$, on which $U_K$ acts as the character $c$, 
where $c$ is one of the $c_i$'s. \\
 
Hence there is a basis $B$ of $E'$ such that $$Mat_B (c^{-1}(t)\tau(t))= \begin{pmatrix}
1 & A_{1,2} (t)     & A_{1,3}(t)     & \dots                      & A_{1,q}  (t)   \\
                    & 1 & A_{2,3} (t)     &  \dots                     & A_{2,q} ( t )   \\
                    &                    & \ddots              &                            & \vdots            \\  
                    &                    &                     & 1      &  A_{q-1,q} ( t ) \\ 
                    &                    &                     &                            & 1                                                                                                    
  \end{pmatrix} $$ for any $t$ in $K^*$, where the $A_{i,j}$'s are smooth functions on $K^*$. So we only have to prove that the $A_{i,j}$'s are polynomials
 of the valuation of $K$ with no constant term.\\
We do this by induction on $q$.\\
 It is obvious when $q=1$.
Suppose the statement holds for $q-1$, and suppose that $E'$ is of dimension $q$, with basis $B=(v_1,\dots,v_{q})$. Considering the two 
$c^{-1}\tau(K^*)$-modules $Vect(v_1,\dots,v_{q-1})$ and  $Vect(v_1,\dots,v_{q})/Vect(v_1)$ of dimension $q-1$, we deduce that for every couple $(i,j)$ different from 
$(1,q)$, there is a polynomial with no constant term $P_{i,j}$ of degree at most $j-i$, such that $A_{i,j}=P_{i,j}\circ v_K$. Now because $c^{-1}\tau$ 
is a representation of $K^*$, and because the $P_{i,j}\circ v_K$'s vanish on $U_K$ for $(i,j)\neq (1,q)$, we deduce that $A_{1,q}$ is a smooth morphism from $(U_K,\times)$ 
to $(\mathbb{C},+)$, which must be zero because $(\mathbb{C},+)$ has no nontrivial compact subgroups. From this we deduce that $A_{1,q}$ is invariant under translation by
 elements of $U_K$ (i.e. $A_{1,q}(\pi_K^k u)= A_{1,q}(\pi_K^k )$ for every $U$ in $U_K$).\\
Denote by $M(k)$ the matrix $Mat_B (c^{-1}\tau(\pi_K^k))$ for $k$ in $\mathbb{Z}$. One has $M(k)=M(1)M(k-1)$ for $k\geq 1$, which in implies 
$A_{1,q}(\pi_K^{k})= \sum_{j=2}^{q-1} P_{1,j}(1)P_{j,q}(k-1) + A_{1,q}(\pi_K^{k-1})+ A_{1,q}(\pi_K)= Q(k)+ A_{1,q}(\pi_K^{k-1})+ A_{1,q}(\pi_K)$ for $Q$
 a polynomial of degree at most $q-2$. This in turn implies that $A_{1,q}(\pi_K^{k})= \sum_{l=1}^{k-1} Q(l) +k A_{1,q}(\pi_K)= R(k)$ for $R$ a polynomial of degree at most $q-1$,
 according to the theory of Bernoulli polynomials, for any $k\geq0$. The same reasoning for $k\leq0$, implies $A_{1,q}(\pi_K^{k})=R'(k)$ for $R'$ a polynomial 
of degree at most $q-1$, for any $k\leq 0$. We need to show that $R=R'$ to conclude.\\ 
We know that $M(k)$ is a matrix whose coefficients are polynomials in $k$ for $k>0$ of degree less that $q-1$, we denote it $P(k)$. The matrix $M(k)$
 has the same property for $k<0$, 
we denote it $P'(k)$.
 Moreover for any $k\geq0$ and $k'\leq 0$, with $k+k'\geq0$, one has $P(k+k')=P(k)P'(k')$. Fix $k >q-1$, then the matrices $P(k +k')$ and $P(k )P'(k')$ 
are equal for $k'$ in $[1-q,0]$, as their coefficients are polynomials in $k'$ with degree at most $q-1$, the equality $P(k +z')=P(k )P'(z')$ holds for any complex 
number $z'$. Now fix such a complex number $z' $, the equality $P(k +z')$ and $P(k )P'(z')$ holds for any integer $k >q-1$, and as both matrices have coefficients which
 are polynomials in $k$, this equality actually holds for any complex number $z$, so that $P(z +z')$ equals $P(z )P'(z')$ for any complex numbers $z$ and $z'$.\\
 As $P(0)=I_q$, we deduce that $P$ and $P'$ are equal on $\mathbb{C}$, and this implies that $R$ is equal to $R'$. 

\end{proof}
 
This proposition has the following consequence:

\begin{prop}\label{zerobehaviour}
 Let $\pi$ be a smooth submodule of finite length of $C^\infty_c(N_n(K)\backslash P_n(K),\psi)$, and let the $c_1,\dots, c_r$ be the central characters of the irreducible
sub factors of $\tau=\pi^{(n-j)}$. Then for any function $W$ in the space of $\pi$, there exist a function $\phi$ in $C_c^\infty(K)$
null at zero, functions $\phi_{k,l}$  in $C_c^\infty(K)$, complex polynomials
$Q_{k,l}$ for $k$ between $1$ and $r$, and $l$ between $1$ and an integer $n_k$, such that one has $W\begin{pmatrix}
t I_j & \\  & I_{n-j}                                                                                                                                           
  \end{pmatrix}= |t|_K^{j(n-j)/2}[\sum_{k=1}^r  c_k(t) [\sum_{l=1}^{n_k} \phi_{k,l}(t)Q_{k,l}(v_K(t))] +\phi(t)]$.
\end{prop}

\begin{proof}
We first remind that the function $\tilde{W}: g\mapsto |g|_K^{-(n-j-1)/2}W \begin{pmatrix}
g & \\  & I_{n-j}                                                                                                                                           
  \end{pmatrix}$ belongs to the space of $\pi_{(n-j-1)}$. We denote by $v$ its image in the space $E$ of $\pi^{(n-j)}$. From Proposition \ref{center-repr},
 the vector $v$ belongs to a finite dimensional $K^*$-submodule $E'$ of $E$, on which $K^*$ acts by a matrix of 
the form determined in Proposition \ref{center-repr}.
We fix a basis $B=(e_{1},\dots,e_{q})$ of $E'$, and denote by $M(a)$ the matrix $M_B(\tau(a))$ (with 
$\tau(a)=\pi^{(n-j)}(aI_j)$), hence we have $\tau(a)e_l= \sum_{k=1}^q M(a)_{k,l}e_k$ for each $l$ between $1$ and $q$.\\
Taking preimages $\tilde{E}_{1},\dots,\tilde{E}_{q}$ of $e_{1},\dots,e_{q}$ in 
$\pi_{(n-j-1)}$, we denote by $\tilde{E}$ the function vector $\begin{pmatrix}
\tilde{E}_{1} \\ \vdots \\  \tilde{E}_{q} \end{pmatrix}$.\\                                                                                                                                       
 As the function 
 $|a|_K^{-j/2}\pi_{(n-j-1)}(aI_j)\tilde{E}_{l}  - \sum_{k=1}^q M(a)_{k,l}\tilde{E}_{k}$ belongs to the kernel of the projection from $\pi_{(n-j-1)}$ to $\pi^{(n-j)}$
(the torsion by $|\ |_K^{-j/2}$ comes from the fact that the projection from the space of $\pi_{(n-j-1)}$ to $\pi^{(n-j)}$ is an intertwining operator of $G_j(K)$-modules,
 if you take the twisted action $|g|_K^{-1/2}\pi_{(n-j-1)}$ on the space of $\pi_{(n-j-1)}$), there is a neighbourhood of zero in $K^j$, such that this function vanishes 
on elements of $G_j(K)$ with last row in this neighbourhood.\\
 In particular, there exists $N_a$ such that the vector function $|a|_K^{-j/2}\pi_{(n-j-1)}(a)\tilde{E}-^t\!M(a)\tilde{E}$ 
vanishes on  $P_K^{N_a}I_j$.\\ 
This implies as in proof of Proposition 2.6. of \cite{CP}, the following claim:

\begin{claim}
There is actually an $N$, such that for every $t$ in $P_K^{N}$, and every $a$ in  $R_K-\{0\}$,
the vector $\tilde{E}(taI_j)$ is equal to ${^t\!M(a)}.|a|_K^{j/2}\tilde{E}(tI_j)$.
\end{claim}

\begin{proof}[Proof of the claim] Indeed, if $U$ is an open compact subgroup of $U_K$, 
such that $\tilde{E}$ and the homomorphism $a \in K^* \mapsto M(a) \in G_q(\mathbb{C})$ are $U$ invariant, and denoting by $u_1,\dots,u_s$ the representatives of $U/U_K$, 
we choose $N$ to be $\max(N_{u_1},\dots,N_{u_q},N_{\pi_K})$. Then for $t$ in $P_K^{N}$, and $a=\pi_K^r u_i u$ in  $R_K$ (with $u$ in $U$), we have 
$\tilde{E}(taI_j)=\tilde{E}(t\pi_K^r u_iI_j)= {^t\!M(u_i)}.|u_i|_K^{j/2}\tilde{E}(t\pi_K^rI_j)$ because $t\pi_K^r$ belongs to $P_K^{N+r}\subset P_K^{N}\subset P_K^{N_{u_i}}$.
But if $r\geq 1$, again one has $\tilde{E}(t\pi_K^rI_j)={^t\!M(\pi_K)}.|\pi_K|_K^{j/2}\tilde{E}(t\pi_K^{r-1}I_j)$ because 
$t \pi_K^{r-1}$ belongs to $P_K^{N+r-1} \subset P_K^{N} \subset P_K^{N_{\pi_K}}$, and repeating this step $r-1$ times, we deduce the equality 
$$\tilde{E}(taI_j)= {^t\!M(u_i)}{^t\!M(\pi_K)}^r |u_i|_K^{j/2} |\pi_K|_K^{rj/2} \tilde{E}(tI_j)={^t\!M(u_i \pi_K^r)}|u_i \pi_K^r|_K^{j/2}\tilde{E}(t)
= {^t\!M(a)} |a|_K^{j/2}\tilde{E}(tI_j).$$ \end{proof}

Denoting by $t_0$ the scalar $\pi_K^N$, the previous claim implies that the vector $\tilde{E}(t)$ is equal to
${^t\!M(t)}|t|_K^{j/2} [{^t\!M(t_0)}^{-1}|t_0|_K^{-j/2}\tilde{E}(t_0 I_j)]$ 
for $t$ in $P_K^N$.
 Let's then choose 
complex numbers $x_1,\dots,x_q$, satisfying $v = (x_1,\dots,x_q)\begin{pmatrix} e_1 \\ \dots \\ e_q \end{pmatrix}$, then there is an integer $N'$ such that we have
$\tilde{W}= (x_1,\dots,x_q)\tilde{E}$ on $P_K^{N'} I_j$, and this in turns implies that if we put $M=max(N,N')$, we have 
$$\tilde{W}(tI_j)=(x_1,\dots,x_q)^t\!M(t)|t|_K^{j/2}V_0 $$ for any $t$ in $P_K^{M}$, and $V_0={^t\!M(t_0)}^{-1}|t_0|_K^{-j/2}\tilde{E}(t_0 I_j)$.\\
We recall that it is a classical fact that because  $\tilde{W}$ is fixed by an open subgroup of $U_j(K)$ ($\pi_{(n-j-1)}$ being smooth) 
and transforms by $\psi$ under left translation by elements of $N_{j}(K)$, that 
the function $\tilde{W}(tI_j)$ vanishes when $t$ is of large absolute value, the preceding equality implies that there is a function $\phi$ in $C_c^\infty(K)$ vanishing
 in a neighbourhood of zero, such that $\tilde{W}(tI_j)$ is equal to $(x_1,\dots,x_q)^t\!M(t)|t|_K^{j/2}V_0 + |t|_K^{j/2}\phi (t)$ 
(as $\phi\mapsto |\ |_K^{j/2}\phi$ is a bijection of the set of functions in $C_c^\infty(K)$ null at zero).\\
 Finally, because the coefficients of $^t\!M(t)$ are of the form a polynomial in $v_K$ multiplied by a central character $c_i$, we deduce that 
there exist functions $\phi_{k,l}$ in $C_c^\infty(K)$, complex polynomials
$Q_{k,l}$ for $k$ between $1$ and $r$, and $l$ between $1$ and an integer $n_k$, such that one has $W\begin{pmatrix}
t I_j & \\  & I_{n-j}                                                                                                                                           
  \end{pmatrix}= |t|_K^{j(n-j)/2}[\sum_{k=1}^r  c_k(t) [\sum_{l=1}^{n_k} \phi_{k,l}(t)Q_{k,l}(v_K(t))] +\phi(t)]$.

\end{proof}

A refinement of the proof of this proposition gives the following:

\begin{thm}\label{restorus}
 Let $\pi$ be a smooth submodule of finite length of $C^\infty_c(N_n(K)\backslash P_n(K),\psi)$, and for $j$ between $1$ and $n-1$, 
let the $c_{1,n-j},\dots, c_{r_j,n-j}$ be the central characters 
of the irreducible sub factors of $\tau=\pi^{(n-j)}$. Then for any function $W$ in the space of $\pi$, the function $$W(t_1,\dots,t_{n-1})=W\begin{pmatrix}
t_1\dots t_{n-1}&  &  &  &  & \\  & t_2 \dots t_{n-1} &  &  &  & \\  &  & \ddots &  & &  \\    &  &  & t_{n-2}t_{n-1} &  & \\ 
   &  & & & t_{n-1} & \\ &  & & &   & 1 \end{pmatrix}$$                                                                                                                                      
 is a linear combination of functions of form $$\prod_{j=1}^{n-1} c_{i_j,j}(t_j)|t_j|^{j(n-j)/2} v_K(t_j)^{m_j}\phi_j(t_j)$$ 
for $i_j$ between $1$ and $r_j$, positive integers $m_j$, and functions $\phi_j$ in $C_c^\infty(K)$.
\end{thm}
\begin{proof} The proof is by induction on $n$.\\
 Let $W$ belong to the space of $\pi$, so that its restriction to $P_n(K)$ belongs to the space of $\pi_{(0)}$, following the beginning of the proof of 
the preceding proposition, we denote by $v$ its image in the space $E$ of $\pi^{(1)}$. Again 
the vector $v$ belongs to a finite dimensional $K^*$-submodule $E'$ of $E$, on which $K^*$ acts by a matrix of 
the form determined in Proposition \ref{center-repr}.
We fix a basis $B=(e_{1},\dots,e_{q})$ of $E'$, and denote by $M(a)$ the matrix $M_B(\tau(a))$ (with 
$\tau(a)=\pi^{(1)}(aI_{n-1})$), hence we have $\tau(a)e_l= \sum_{k=1}^q M(a)_{k,l}e_k$ for each $l$ between $1$ and $q$.\\
Taking preimages $\tilde{E}_{1},\dots,\tilde{E}_{q}$ of $e_{1},\dots,e_{q}$ in 
$\pi_{(0)}$, we denote by $\tilde{E}$ the function vector $\begin{pmatrix}
\tilde{E}_{1} \\ \vdots \\  \tilde{E}_{q} \end{pmatrix}$.\\                                                                                                                                       
There is a neighbourhood of zero in $K^{n-1}$, such that the function $|a|_K^{-(n-1)/2}\pi_{(0)}(a)\tilde{E}_{l}  - \sum_{k=1}^q M(a)_{k,l}\tilde{E}_{k}$  vanishes 
on elements of $G_{n-1}(K)$ with last row in this neighbourhood.
In particular, there exists $N_a$ such that for every $(t_1,\dots,t_{n-2})$ in $(K^*)^{n-2}$, the vector function 
$$|a|_K^{-(n-1)/2}\pi_{(0)}(aI_{n-1})\tilde{E}(t_1,\dots,t_{n-1}) -^t\!M(a)\tilde{E}(t_1,\dots,t_{n-1})$$ 
vanishes when $t_{n-1}$ belongs to $P_K^{N_a}I_j$ (here by $\tilde{E}(t_1,\dots,t_{n-1})$, we mean $$\tilde{E}\begin{pmatrix}
t_1\dots t_{n-1}&  &  &  &  & \\  & t_2 \dots t_{n-1} &  &  &  & \\  &  & \ddots &  & &  \\    &  &  & t_{n-2}t_{n-1} &  & \\ 
   &  & & & t_{n-1} & \\ &  & & &   & 1 \end{pmatrix},$$ hence $\pi_{(0)}(aI_{n-1})\tilde{E}(t_1,\dots,t_{n-1})=\tilde{E}(t_1,\dots,t_{n-1}a)$).\\ 
As in the proof of Proposition \ref{zerobehaviour}, we deduce that there is an integer $N_1$, and an element $t_0$ of $K^*$, 
such that for every $(t_1,\dots,t_{n-2})$ in $(K^*)^{n-2}$, the vector $\tilde{E}(t_1,\dots,t_{n-1})$ is equal to 
$^t\!M(t_{n-1})|t_{n-1}|_K^{(n-1)/2}[{^t\!M(t_{0}^{-1})}|t_{0}|_K^{-(n-1)/2}]\tilde{E}(t_1,\dots, t_{n-2},t_{0})$ for any $t_{n-1}$ in $P_K^{N_1}$.\\
If the image $v$ of $W$ in $\pi^{(1)}$ is equal to $x_1 e_1+\dots+x_q e_q$, there is also an integer $N_2$, 
such that for every $(t_1,\dots,t_{n-2})$ in $(K^*)^{n-2}$, the function $W(t_1,\dots,t_{n-1})-(x_1,\dots,x_q)\tilde{E}(t_1,\dots,t_{n-1})$ vanishes when $t_{n-1}$ belongs
to $P_K^{N_2}$. Hence taking $N=(max(N_1,N_2)$, for any $t_{n-1}$ in $P_K^{N}$, and any $(t_1,\dots,t_{n-2})$ in $(K^*)^{n-2}$ the function $W(t_1,\dots,t_{n-1})$ 
is equal to $$(x_1,\dots,x_q) ^t\!M(t_{n-1})|t_{n-1}|_K^{(n-1)/2}[{^t\!M(t_{0}^{-1})}|t_{0}|_K^{-(n-1)/2}]\tilde{E}(t_1,\dots, t_{n-2},t_{0}).$$ \\
But the functions $|\ |_K^{-k/2}\pi_{(0)}(t_{0}I_{n-1})\tilde{E}_i$ belong to the smooth submodule of finite length $\Phi^-(\pi)$ of $C^\infty_c(N_{n-1}(K)\backslash P_{n-1}(K),\psi)$,
 so by induction hypothesis, the functions $\tilde{E}(t_1,\dots, t_{n-2},t_{0})$ are sums of functions of the form 
$\prod_{j=1}^{n-2} c_{i_j,j}(t_j)|t_j|^{j(n-j)/2} v_K(t_j)^{m_j}\phi'_j(t_j)$
for $i_j$ between $1$ and $r_j$, positive integers $m_j$, and functions $\phi'_j$ in $C_c^\infty(K)$. This in turn implies that the function 
$$S(t_1,\dots,t_{n-1})=(x_1,\dots,x_q) ^t\!M(t_{n-1})|t_{n-1}|_K^{(n-1)/2}[{^t\!M(t_{0}^{-1})}|t_{0}|_K^{-(n-1)/2}]\tilde{E}(t_1,\dots, t_{n-2},t_{0})$$ is also
 a sum of functions of the form 
$\prod_{j=1}^{n-1} c_{i_j,j}(t_j)|t_j|^{j(n-j)/2} v_K(t_j)^{m_j}\phi_j(t_j)$ for $i_j$ between $1$ and $r_j$, positive integers $m_j$, and functions 
$\phi_j$ in $C_c^\infty(K)$.\\
Reminding that there is an integer $N'$, such that for any $(t_1,\dots,t_{n-2})$ in $(K^*)^{n-2}$, and any $t_{n-1}$ of absolute value greater than $q_K^{N'}$, 
both $W(t_1,\dots,t_{n-1})$ and $S(t_1,\dots,t_{n-1})$ are zero, we deduce that the difference of the two
functions is a smooth function $\phi(t_1,\dots,t_{n-1})$ on $(K^*)^{n-1}$ which vanishes when $t_{n-1}$ has absolute value outside $[q_K^{-N},q_K^{N'}]$.
 Moreover there is a compact subgroup $U$ of $U_K$ independent of $(t_1,\dots,t_{n-1})$ such that both functions (hence $\phi$) are invariant when $t_{n-1}$ 
is multiplied by an element of $U$. Denoting by $(t_\alpha)_{\alpha \ \in \ A}$ a finite set of representatives of $$\{ t\ | q_K^{-N}\leq |t_k|_K \leq q_K^{N'}\}/U,$$
 this implies that $\phi(t_1,\dots,t_{n-1})$ is equal to $\sum_{\alpha \ \in \ A} \phi(t_1,\dots,t_{n-2},t_{\alpha})\mathbf{1}_{t_\alpha U}(t_{n-1})$,
 which we can always write as $\sum_{\alpha \ \in \ A} \phi(t_1,\dots,t_{n-2},t_{\alpha})|t_{n-1}|^{(n-1)/2}\phi_\alpha(t_{n-1})$ with 
$\phi_\alpha=|\ |_K^{-(n-1)/2}\mathbf{1}_{t_\alpha U}$ in $C_c^\infty(K)$ and null at zero.\\ 
Finally as each function $ \phi(t_1,\dots,t_{n-2},t_{\alpha})$ is equal to $W(t_1,\dots,t_{\alpha})-S(t_1,\dots,t_{\alpha})$, by induction hypothesis again, it is a sum 
of functions of the form $\prod_{j=1}^{n-2} c_{i_j,j}(t_j)|t_j|^{j(n-j)/2} v_K(t_j)^{m_j}\phi''_j(t_j)$
for $i_j$ between $1$ and $r_j$, positive integers $m_j$, and functions $\phi''_j$ in $C_c^\infty(K)$, and the statement of our proposition follows. 
\end{proof}

We are now able to prove Proposition \ref{l2}.\\

Indeed, let $W$ belong to the space of $\pi$ as in the statement of Proposition \ref{l2}, first we notice the equality
$$\int_{N_n(K)\backslash P_n(K)} |W(p)|^2 dp=\int_{N_{n-1}(K)\backslash G_{n-1}(K)} |W(g)|^2 dg.$$

Now the Iwasawa decomposition reduces the convergence of this integral to that of $$\int_{A_{n-1}(K)} |W(a)|^2 \Delta_{B_{n-1}(K)}(a) d^*a.$$
Using coordinates $(t_1,\dots,t_{n-1})= \begin{pmatrix}
t_1\dots t_{n-1}&  &  &  &  & \\  & t_2 \dots t_{n-1} &  &  &  & \\  &  & \ddots &  & &  \\    &  &  & t_{n-2}t_{n-1} &  & \\ 
   &  & & & t_{n-1} & \\ &  & & &   & 1 \end{pmatrix}$ of ${A_{n-1}(K)}$, the function $\Delta_{B_{n-1}(K)}(t_1,\dots,t_{n-1})$ is equal to 
$\prod_{j=1}^{n-1} |t_j|_K^{-j(n-j-1)}$.\\
According to Proposition \ref{restorus} the function $|W(t_1,\dots,t_{n-1})|^2$ is majorized by a sum of functions of the form 
$$\prod_{j=1}^{n-1} |c_{k_j,j}|(t_j)|c_{l_j,j}|(t_j)|t_j|_K^{j(n-j)} v_K(t_j)^{m_j}\phi_j(t_j)$$ for $c_{k_j,j}$ and $c_{l_j,j}$ central characters 
of irreducible sub factors of $\pi^{(n-j)}$, $m_j$ non negative integers, and $\phi_j$ non negative functions in $C_c^\infty(K)$.\\
Hence our integral will converge if so do the integrals 
$$\int_{A_{n-1}(K)} \prod_{j=1}^{n-1} |c_{k_j,j}|(t_j)|c_{l_j,j}|(t_j)|t_j|_K^{j} v_K(t_j)^{m_j}\phi_j(t_j)d^*t_1\dots d^*t_n,$$ i.e. if so does the integral 
 $\int_{K^*} |c_{k_j,j}|(t_j)|c_{l_j,j}|(t_j)|t_j|_K^{j} v_K(t_j)^{m_j}\phi_j(t_j)d^*t_j$ for any $j$ between $1$ and $n-1$.\\
But our assertion on the central exponents of the shifted derivatives, insures that we have $|c_{k_j,j}|(t_j)=|t_j|_K^{r_1}$ and $|c_{l_j,j}|(t_j)=|t_j|_K^{r_2}$,
 with $r_1$ and $r_2$ both strictly greater than $-j/2$, and it is than classical that the integral 
$\int_{K^*} |t_j|_K^{j-(r_1+r_2)} v_K(t_j)^{m_j}\phi_j(t_j)d^*t_j$ converges.\\
This proves Proposition \ref{l2}.

\end{document}